\documentclass[fleqn,12pt,twoside]{article}
\usepackage[headings]{espcrc1}
\usepackage{amsmath}
\usepackage{graphicx}
\PassOptionsToPackage{ctagsplt,righttag}{amsmath}
\usepackage{multirow}

\newcommand{\beq}{\begin{equation}}
\newcommand{\eeq}{\end{equation}}
\newcommand{\bea}{\begin{eqnarray}}
\newcommand{\eea}{\end{eqnarray}}

\def\pt{\partial t}
\def\px{\partial x}
\def\py{\partial y}

\def\tu{\tilde{u}}
\def\he{\hat{e}}
\def\dt{\Delta t}
\def\dx{\Delta x}
\def\dy{\Delta y}
\def\dto{\dt_{opt}}

\def\eg{\emph{e.g. }}
\def\ie{\emph{i.e. }}
\def\etal{\emph{et al.}}

\begin{document}
\bibliographystyle{unsrt}

\title{Boosting the Accuracy of Finite Difference Schemes via
       Optimal Time Step Selection and Non-Iterative Defect Correction}

\author{Kevin T. Chu\address{Institute of High Performance Computing, A*STAR, Singapore, Singapore}$^,$\address{Vitamin D, Inc., Menlo Park, CA 94025, United States}
}

\runtitle{Optimal Time Step Selection and Non-Iterative Defect Correction for 
          FD Schemes}
\runauthor{K.T. Chu}

\maketitle

\noindent \rule{6.3in}{1pt}

\begin{abstract}
In this article, we present a simple technique for boosting the order of 
accuracy of finite difference schemes for time dependent 
partial differential equations by optimally selecting the time step used
to advance the numerical solution and adding defect correction terms in a
non-iterative manner.  The power of the technique is its ability to extract 
as much accuracy as possible from existing finite difference schemes with 
minimal additional effort.  Through straightforward numerical analysis 
arguments, we explain the origin of the boost in accuracy and estimate the 
computational cost of the resulting numerical method.  We demonstrate the 
utility of optimal time step (OTS) selection combined with non-iterative 
defect correction (NIDC) on several different types of finite difference
schemes for a wide array of classical linear and semilinear PDEs in one and 
more space dimensions on both regular and irregular domains.  
\end{abstract}

\noindent \rule{6.3in}{1pt}

% KEYWORDS
% optimal time step; non-iterative defect correction; finite difference 
% schemes; high-order accurate numerical methods; time dependent PDEs

\section{Introduction}
High-order numerical methods for partial differential equations (PDEs) will 
always be valuable for increasing the computational efficiency of numerical 
simulations.  Thus, it is not at all surprising that a great deal of effort in 
numerical PDEs continues to be focused on the development of high-order 
numerical 
schemes~\cite{spotz_2001,bruger_2005,gibou_2005,ito_2005,shukla_2005,heidenreich_2007,shukla_2007}.  
Typically, high-order accuracy is achieved by constructing
schemes that are formally high-order accurate.  However, high-order 
accuracy can also be obtained by using formally low-order schemes in clever 
ways (\eg Richardson extrapolation and compact finite difference 
schemes~\cite{spotz_2001}).  When possible, the latter approach can be a 
powerful way to boost the accuracy of a numerical method \emph{without} 
introducing too much additional algorithmic (and programming) complexity.

Optimization of numerical parameters (\eg $\dx$ and $\dt$) is an important, 
but relatively uncommon, approach for maximizing the accuracy of numerical 
schemes.  Often, all combinations of numerical parameters are considered 
equally good.  For some problems, however, optimizing the numerical parameters 
can significantly increase the accuracy of the numerical solution.  For 
example, the well-known unit CFL condition, $\dt = \dx/a$, for the first-order 
upwind forward Euler discretization of the linear advection equation 
completely eliminates the numerical error introduced by the 
scheme~\cite{leveque_book_2002}.  Another well-known example is the boost in 
the order of accuracy that results from setting the time step to 
$\dt = \dx^2/6D$ when the 1D diffusion equation is solved using forward Euler 
time integration with the standard second-order central difference 
approximation for the Laplacian.  As a final example, Zhao showed that an 
appropriate choice of the weight factor $\theta$ in the theta-method can 
significantly decrease the error in finite element solutions to the diffusion 
equation~\cite{zhao_2006}. 

In this article, we present a simple technique for boosting the order of 
accuracy of finite difference schemes for time dependent PDEs by optimally 
choosing the time step and, when necessary, adding defect correction terms
in a \emph{non-iterative} manner.  This technique is a systematic 
generalization of the ideas underlying some of the examples mentioned above.
Optimal time step (OTS) selection is based on the observation that a carefully 
chosen time step can simultaneously eliminate low-order terms in both the
spatial and temporal the discretization errors.  It combines the well-known 
approach of using knowledge of the structure of the PDE to reduce numerical 
error with the less-common notion that an optimal choice of numerical 
parameters can boost the effectiveness of algorithmically simple and 
computationally low-cost numerical schemes.  For some problems, OTS alone is 
sufficient to allow a formally low-order method to deliver high-order 
accuracy.  However, addition of defect correction terms is generally required 
to eliminate residual terms in the leading-order truncation error.  It is 
important to emphasize that, unlike traditional defect/deferred correction 
methods~\cite{pereyra_1968,stetter_1978,gustafsson_2002,kress_2002,kress_2006} 
which use defect correction terms to iteratively refine the solution at 
each time step, we take a non-iterative defect correction (NIDC) approach
that incorporates defect correction terms directly into the original finite 
difference scheme.

As a technique for enhancing numerical methods, OTS-NIDC has several desirable 
features.  Like any technique that leads to high-order methods, OTS-NIDC
produces schemes that greatly reduce the computational cost required to obtain 
a numerical solution of a desired accuracy.  However, its real power comes 
from the fact that it achieves high-order accuracy while being based solely on 
very simple, formally low-order finite difference schemes.  In other words, 
OTS-NIDC allows us to extract as much accuracy as possible from an existing 
numerical scheme with minimal additional work.  For instance, we will see how 
OTS-NIDC makes it possible to achieve 4th-order accuracy for the diffusion 
equation in any number of space dimensions using only a second-order stencil 
for the Laplacian and explicit time integration.  
Moreover, high-order accuracy is not restricted to problems 
on simple, rectangular domains; irregular domains are straightforward to 
handle by appropriately setting the values in ghost 
cells~\cite{gibou_2005,ito_2005,fedkiw_1999,osher_fedkiw_book}.  
OTS-NIDC is even beneficial when using finite difference schemes to solve some 
nonlinear PDEs.  While there are certainly limitations to OTS-NIDC, we will 
see that there are many problems where it is useful.  

This article is organized as follows.  In the remainder of this section,
we compare OTS-NIDC to a few related numerical techniques.
In Section~\ref{sec:OTS-NIDC}, we present the theoretical foundations of 
OTS-NIDC and show how it boosts the order of accuracy for formally low-order 
numerical schemes.  
In Section~\ref{sec:applications_1d}, we demonstrate the utility of 
OTS-NIDC to problems in one space dimension by applying it to finite 
difference schemes for several classical linear and nonlinear PDEs.  
In Section~\ref{sec:applications_multidim}, we show that 
OTS-NIDC is useful for problems in multiple space dimensions by applying it
the 2D linear advection equation and the 2D diffusion equation.  To illustrate 
the ease with which high-order accuracy can be achieved for problems on 
arbitrary domains, the 2D diffusion equation is solved on both regular 
and irregular domains.  
For each PDE in Sections~\ref{sec:applications_1d} and 
\ref{sec:applications_multidim}, we provide the numerical analysis required 
to calculate the optimal time step size, the defect correction terms, and the 
order of accuracy.  Finally, we offer some concluding remarks in 
Section~\ref{sec:summary}.

\subsection{Comparison with Traditional Defect/Deferred Correction Methods}
Defect/deferred correction is a well-known strategy for boosting the order of 
accuracy of finite difference schemes for time dependent PDEs.  It usually 
takes one of two forms: 
(1) iterative refinement of the numerical solution~\cite{pereyra_1968,stetter_1978,gustafsson_2002,kress_2002,kress_2006}
or 
(2) derivation of high-order compact schemes for the spatial derivative 
operator treating time derivatives as source 
terms~\cite{spotz_2001,ito_2005,heidenreich_2007}. 

For iterated defect/deferred correction, the numerical solution at each 
time step is computed by iteratively adding correction terms that are 
calculated from lower accuracy iterates.  When used to improve the 
temporal order of accuracy, the iterative approach can help to ensure 
stability of the time integration scheme~\cite{kress_2006}.  Although there 
are instances where iterative defect correction has been used to eliminate 
both spatial and temporal errors~\cite{gustafsson_2002}, it is usually used to 
improve only one of these errors~\cite{pereyra_1968,kress_2002,kress_2006}.
Note that unlike iterated defect correction, OTS-NIDC requires no iteration 
during each time step because it computes correction terms using only the 
solution at the previous time step(s) (as opposed to requiring an estimate 
of the solution at the current time).

High-order compact schemes focus on eliminating spatial discretization errors 
when deriving the semi-discrete equations.  When applied to time dependent 
PDEs, a high-order compact scheme for the spatial derivative operator is 
derived treating time derivatives as source 
terms~\cite{spotz_2001,ito_2005,heidenreich_2007}.  
After deriving the high-order spatial discretization, any time differencing 
scheme may be applied to discretize the time derivatives.  While this 
approach is very general, it really only provides a mechanism for reducing 
spatial discretization error.  Temporal discretization errors are completely 
controlled by the choice of time differencing scheme.  Note that because 
high-order compact schemes for the spatial derivative operators involve spatial 
derivatives of the source term, the resulting schemes are always effectively 
implicit (even if all of the spatial derivative terms are treated explicitly).  

OTS-NIDC has several advantages over both of these defect correction 
techniques.  First, OTS-NIDC is generally more computationally efficient than 
both iterated defect correction and the high-order compact scheme approach
because it requires less memory and, for many problems, admits simple 
explicit time integration schemes.  Second, OTS-NIDC inherently controls both 
spatial and temporal numerical errors at the same time.  Finally, OTS-NIDC 
does not rely solely on the addition of finite differences of solution values 
to cancel out discretization errors.  Instead, it takes advantage of 
relationships between spatial derivatives, temporal derivatives and source 
terms to eliminate the need for finite differences approximations of 
high-order derivatives of the solution.  Compared with iterated defect 
correction, OTS-NIDC makes it easier to avoid instabilities that arise from 
the presence of finite difference discretizations of high-order temporal 
derivatives and may eliminate some of the correction terms.  Comparing to the
high-order compact scheme approach, the optimization of numerical parameters 
in OTS-NIDC means that mixed partial derivatives do not need to appear in the 
fully discrete equations.

\subsection{Comparison with Adaptive Time Stepping Methods}
OTS selection is distinct and separate from traditional adaptive time stepping 
techniques~\cite{iserles_book,shampine_2005} that are used in the 
context of the method of lines~\cite{iserles_book,gko_book}.  
While both methods reduce numerical error by controlling the time step, the 
numerical errors they control are fundamentally different.  Adaptive time 
stepping can only be used to reduce temporal discretization errors because it 
controls the errors that arise when solving the coupled system of ODEs for the 
semi-discrete approximation to the PDE; the spatial discretization errors 
are completely controlled by the finite difference scheme selected to 
approximate the spatial derivatives.  
In contrast, optimal time step selection simultaneously reduces 
\emph{both} the spatial and temporal discretization errors because it
uses information from the PDE to choose the time step.  Another important 
distinction between the two techniques is that optimal time step selection 
uses a \emph{fixed} time step which makes it significantly less 
computationally complex than adaptive time stepping methods.

\section{\label{sec:OTS-NIDC}
         Theoretical Foundations of OTS Selection and NIDC}
OTS-NIDC is not by itself a method for constructing finite difference 
schemes.  Rather, it is a way to enhance to the performance of existing 
finite difference schemes by carefully choosing the time step and adding
defect correction terms to eliminate low-order numerical errors.  
The two key ideas underlying OTS-NIDC are 
(1) a judicious choice of time step can be used to eliminate the leading-order 
(and possibly higher-order) terms in the error
and
(2) the PDE provides valuable insight into the relationships between terms
in the discretization error. 
In this section, we present and analyze the technique of OTS-NIDC for finite 
difference schemes constructed to solve scalar evolution equations of the form 
\beq
  \frac{\partial u}{\partial t} = 
    F\left(x, u(x), D u(x), D^2 u(x), \ldots \right), 
\eeq
where the right hand side is an arbitrary function of the solution $u$, its 
spatial derivatives, and the independent variable $x$.
Because this class of time dependent PDEs is so prevalent and important, 
restricting our attention to problems of this form is not a serious 
limitation.  The same principles apply to finite difference schemes for 
general PDEs.

\subsection{\label{sec:ots_model_pde} 
            Optimal Time Step for a Model PDE}
Let us begin by considering a finite difference approximation for the
following time dependent PDE in one spatial dimension: 
\beq
  \frac{\partial u}{\partial t} = A \frac{\partial^n u}{\partial x^n},
  \label{eq:model_PDE}
\eeq
where $A$ is a constant of the appropriate sign so that~(\ref{eq:model_PDE}) 
is well-posed.  Because we will be boosting the accuracy of the finite 
difference approximation by optimizing the time step, there is no reason to 
explicitly construct the scheme to be high-order.  In fact, using a high-order 
finite difference scheme only complicates the analysis required to calculate 
the optimal time step.  

To compute the optimal time step, we first analyze the discretization errors 
(both spatial and temporal) for the finite difference approximation, 
explicitly keeping track of the leading-order terms.  In contrast to 
conventional error analysis, we do not sweep all of the errors under big-$O$
notation.  Rather, we take advantage of the fact that derivation of the 
leading-order terms in the error is straightforward for finite difference
schemes through the use of Taylor series expansions.  Because each term in 
the error is proportional to a partial derivative of the solution, we can 
use the PDE to relate the leading-order terms to each other.  
For many finite difference schemes, this careful analysis yields a direct 
relationship between the leading-order spatial and temporal errors that allows 
them to be combined as a single term of the form:
\beq
  (L u) P(\dx, \dt) ,
  \label{eq:leading_order_error_model_PDE_general}
\eeq
where $L$ is a differential operator and $P$ is a simple function of the 
numerical parameters $\dx$ and $\dt$.  By selecting the numerical parameters 
so that $P = 0$, we can completely eliminate the leading-order term in the 
discretization error.  Because time dependent PDEs are often solved by 
stepping in time, it is natural to choose the time step to be a function of 
the grid spacing.  The time step selected in this way is the
\emph{optimal time step} for the finite difference scheme.  Using the optimal
time step for time stepping yields a numerical solution with an order 
of accuracy that is \emph{higher} than the formal order of the original finite 
difference scheme.

An important class of finite difference schemes for time dependent PDEs are 
those based on the method of lines.   When the method of lines approach is 
used for (\ref{eq:model_PDE}), $P$ often takes a particularly simple 
form:
\beq
  P(\dx, \dt) = (\alpha \dx^r - \beta \dt^s) \dt,
  \label{eq:leading_order_error_model_PDE_MOL}
\eeq
where $r$ and $(s+1)$ are the orders of the leading spatial and temporal 
discretization errors, and $\alpha$ and $\beta$ are constants that 
depend on the details of the PDE and finite difference scheme.  Setting 
$P = 0$, we find that the optimal time step is given by\footnote{The $\dt = 0$ 
root of $P$ is discarded because it is not numerically meaningful.}
\beq
  \dto = \left(\alpha/\beta \right)^{1/s} \dx^{r/s}.
  \label{eq:optimal_time_step}
\eeq
Notice that the optimal time step has a similar functional dependence on the 
grid spacing as a time step that satisfies a stability constraint.  

An important assumption in the above analysis is that $\alpha$ and $\beta$ 
have the same sign.  If this condition is not satisfied, (\eg finite
difference scheme for the diffusion equation based on the standard central
difference stencil for the Laplacian and backward Euler time integration),
then OTS selection is not applicable because none of the zeros of 
(\ref{eq:leading_order_error_model_PDE_MOL}) is positive.  Fortunately, 
in this situation, it is usually possible to find a closely related finite 
difference scheme which \emph{does} have a positive optimal time step.

\subsubsection*{\label{sec:error_analysis} 
            Error Analysis}
The form of the local error in 
(\ref{eq:leading_order_error_model_PDE_MOL}) follows from the fact 
that the method of lines first approximates the PDE as a system of ODEs in 
time by only discretizing the spatial derivatives: 
\beq
{\mathbf u}_t = F({\mathbf u}),
\label{eq:method_of_lines}
\eeq
where ${\mathbf u}$ is the vector of values of $u$ at the grid points and
$F({\mathbf u})$ is the finite difference approximation of the spatial 
derivative operator acting on ${\mathbf u}$.  The second term in the local 
error (\ref{eq:leading_order_error_model_PDE_MOL}) arises directly from 
temporal discretization errors associated with the numerical scheme used to 
integrate the system of ODEs in time.  The first term arises because $F$ is 
only an approximation of the spatial derivative operators in the PDE.  As a 
result, the system of ODEs~(\ref{eq:method_of_lines}) has an error due to 
the spatial discretization at \emph{every} point in time.  For a spatial 
discretization error of order $r$, the error in (\ref{eq:method_of_lines}) is 
$O(\dx^r)$, so the error accumulated during each time step due to 
spatial discretization error is $O(\dt \dx^r)$.  Because the method of 
lines is so frequently used to construct finite difference approximations for 
time dependent PDEs, we will focus primarily on finite difference schemes of 
this type for the remainder of the article.  It is straightforward to extend 
the results we obtain to more general finite difference schemes (\eg the
Lax-Wendroff scheme for the linear advection equation). 

To analyze the accuracy of the finite difference scheme 
for~(\ref{eq:model_PDE}) when the optimal time step is used, we examine the 
higher-order terms in the discretization error.  Elimination of the leading 
terms in the discretization error leaves a local error of the form 
$O(\dt \dx^p) + O(\dt^{q+1})$, where $p>r$ and $(q+1) > (s+1)$ are the orders 
of the spatial and temporal discretization errors, respectively, after the 
leading-order terms in the error have been eliminated.  
From the local error, we can compute the global error of the numerical 
solution over extended intervals of time by using the heuristic 
argument that the global error is equal to the local error divided by 
$\dt$~\cite{gko_book}:
\beq
e = O(\dx^p) + O(\dt^q).
\label{eq:global_error_ots}
\eeq
Because $p > r$ and $q > s$, we see that using an optimal time step 
boosts the overall accuracy of the finite difference scheme.

It is important to emphasize that $\dx$ and $\dt$ are \emph{not} independent 
in~(\ref{eq:global_error_ots}) because the time step must be set to $\dto$ 
in order to derive the error estimate.  In other words, there is really only 
one parameter, $\dx$, that controls the numerical accuracy of the finite 
difference scheme (see Appendix~\ref{app:formal_vs_practical_accuracy} for 
further discussion about this important point).  Since $\dto = O(\dx^{r/s})$, 
the global error (\ref{eq:global_error_ots}) reduces to
\beq
e = O \left( \dx^{\min(p,rq/s)} \right).
\label{eq:global_error_ots_simplified}
\eeq

\subsubsection*{Stability Considerations}
We have not yet considered the issue of numerical stability of finite
difference schemes when the time step is set to $\dto$.  For explicit time 
integration schemes, OTS selection is ineffective when the optimal time step 
is greater than the smallest stable time step.  In this situation, it is 
necessary to replace the finite difference approximation of the PDE.  
For example, if the diffusion equation is solved using a 
finite difference scheme based on forward Euler time integration and a 
discretization of the Laplacian with an isotropic leading-order truncation 
error (see Section~\ref{sec:diffusion_eqn_2d}), the optimal time step becomes 
unstable when the dimension of space is greater than 3.  To make OTS selection 
useful for diffusion problems in 4 and more space dimensions, we can switch 
to the DuFort-Frankel scheme, which is explicit and unconditionally 
stable~\cite{gko_book}.

\subsubsection*{\label{sec:computational_performance} 
                Computational Performance}
For many problems, the optimal time step is on the order of the maximum 
time step allowed by the stability constraint for an explicit method.  
As a result, using the optimal time step to advance the finite
difference scheme in time may appear to be overly restrictive.
However, this apparent drawback is more than 
compensated by the boost in the order of accuracy.  For example, when solving 
the diffusion equation using a forward Euler time integration scheme with a 
second-order central difference stencil for the Laplacian, we obtain a 
fourth-order scheme when $\dt$ is set equal to the optimal time step 
(discussed further in Sections~\ref{sec:diffusion_eqn_1d} and
\ref{sec:diffusion_eqn_2d}).  
As Table~\ref{tab:comp_perf_vs_err} shows, this boost in accuracy leads to a
scheme that is computationally cheaper than traditional methods that use 
the same low-order finite difference stencils.  The performance gain for 
problems in higher spatial dimensions is even more impressive.  As we can see 
in Table~\ref{tab:comp_perf_vs_dim}, the gap between the performance of a
simple forward Euler scheme with the optimal time step and the Crank-Nicholson
method grows with the spatial dimension of the problem.  

\begin{table}[tb]
\caption{
Computational cost of various finite difference schemes for the 1D 
diffusion equation as a function of the global numerical error $e$.
For all of the schemes, the standard second-order central difference 
stencil is used to discretize the Laplacian.  
Note that the time step for the backward Euler scheme is chosen so that the
temporal error is subdominant to the spatial error and that the computation 
time for the Crank-Nicholson method assumes an $O(N) = O \left( 1/\dx\right)$ 
matrix inversion algorithm. 
}
\label{tab:comp_perf_vs_err} 
\renewcommand{\arraystretch}{1.5}
\centering
\begin{tabular}{lcccc}
  \hline
  {\bf Numerical Scheme} & $\dx$ 
  & $\dt$
  & {\bf Memory}
  & {\bf Compute Time}
  \\
  \hline 
  Forward Euler    & $O\left( e^{1/2} \right)$ 
                   & $O\left( e \right)$ 
                   & $O\left( e^{-1/2} \right)$ 
                   & $O\left( e^{-3/2} \right)$ \\
  Backward Euler   & $O\left( e^{1/2} \right)$ 
                   & $O\left( e \right)$ 
                   & $O\left( e^{-1/2} \right)$ 
                   & $O\left( e^{-3/2} \right)$ \\
  Crank-Nicholson  & $O\left( e^{1/2} \right)$ 
                   & $O\left( e^{1/2} \right)$ 
                   & $O\left( e^{-1/2} \right)$ 
                   & $O\left( e^{-1} \right)$ \\
  Forward Euler with OTS-NIDC  & $O\left( e^{1/4} \right)$ 
                   & $O\left( e^{1/2} \right)$ 
                   & $O\left( e^{-1/4} \right)$ 
                   & $O\left( e^{-3/4} \right)$ \\ 
  \hline
\end{tabular}
\end{table}

\begin{table}[tb]
\caption{
Computational cost of the forward Euler scheme with OTS-NIDC and
the Crank-Nicholson method for solving the diffusion equation as a function 
of the global numerical error $e$.
Note that the computation time for the Crank-Nicholson method assumes an 
$O(N) = O \left( 1/\dx^d \right)$ matrix inversion algorithm. 
A comparison of computational costs in higher spatial dimensions continues 
along the same trends but requires replacement of the forward Euler time 
integration with an alternative for which the optimal time step is stable 
(\eg the DuFort-Frankel scheme). 
}
\label{tab:comp_perf_vs_dim}
\renewcommand{\arraystretch}{1.5}
\centering
\begin{tabular}{ccccc}
  \hline
  & \multicolumn{2}{c}{\bf Forward Euler with OTS-NIDC} 
  & \multicolumn{2}{c}{\bf Crank-Nicholson} \\
  \cline{2-3} \cline{4-5} 
    {\bf Dimensions} & {\bf Memory} & {\bf Compute Time} 
  & {\bf Memory} & {\bf Compute Time} \\
  \hline 
  $1$ & $O\left( e^{-1/4} \right)$ 
      & $O\left( e^{-3/4} \right)$ 
      & $O\left( e^{-1/2} \right)$ 
      & $O\left( e^{-1} \right)$ \\ 
  $2$ & $O\left( e^{-1/2} \right)$ 
      & $O\left( e^{-1} \right)$ 
      & $O\left( e^{-1} \right)$ 
      & $O\left( e^{-3/2} \right)$ \\ 
  $3$ & $O\left( e^{-3/4} \right)$ 
      & $O\left( e^{-5/4} \right)$ 
      & $O\left( e^{-3/2} \right)$ 
      & $O\left( e^{-2} \right)$ \\
  \hline 
\end{tabular}
\end{table}

\subsection{\label{sec:ots_nidc_general_1d_pdes} 
            OTS-NIDC for General PDEs in One Space Dimension} 
In the previous section, we analyzed optimal time step selection in the 
context of a simple model PDE in one space dimension.  For more general PDEs, 
OTS selection continues to be useful, but it needs to be augmented with 
defect correction (incorporated in a non-iterative manner) to achieve 
the same improvement in the order of accuracy.  In this section, we discuss 
these correction terms in the context of general constant coefficient, linear 
PDEs and semilinear PDEs.

\subsubsection*{\label{sec:ots_linear_pde} 
            General Constant Coefficient, Linear PDEs} 
Let us consider a general constant coefficient, linear PDE of the form
\beq
  \frac{\partial u}{\partial t} = 
  \sum_{k=1}^n A_k \frac{\partial^k u}{\partial x^k} + f(x,t)
  \label{eq:linear_PDE}.
\eeq
Even if we use the PDE to relate terms in the discretization error, it is
unlikely that all of the dominant terms in the error can be combined into 
a single term.   As a result, the leading-order discretization error is
generally a sum of several terms:
\beq
  \dt\left[ (L u) (\alpha \dx^r - \beta \dt^s) 
  + \sum_k (L_k u) \dx^r 
  + \sum_k (L'_k u) \dt^s 
  + (G f) \dt^s \right]
  \label{eq:leading_order_error_linear_PDE_MOL}
\eeq
where $L$, $u$, $\alpha$, $\beta,$ $r$, and $s$ are defined as in
(\ref{eq:leading_order_error_model_PDE_general}) and
(\ref{eq:leading_order_error_model_PDE_MOL}), 
$L_k$ and $L'_k$ are the linear spatial differential operators associated 
with the leading-order terms in the discretization error that are not involved 
in the choice of the optimal time step, and $G$ is a spatio-temporal 
differential operator that acts on the source term $f$.  
Note that neither $L_k$ nor $L'_k$ involve any temporal differential 
operators because any that arise can be eliminated using the 
PDE~(\ref{eq:linear_PDE}).  

From~(\ref{eq:leading_order_error_linear_PDE_MOL}), it is clear that setting 
$\dt = \dto$ is not enough to completely eliminate the 
leading-order error.  To eliminate the vestigial errors, we need to add 
defect correction terms to cancel out the second, third, and fourth terms in 
(\ref{eq:leading_order_error_linear_PDE_MOL}).  At each time step, the second 
and third terms in (\ref{eq:leading_order_error_linear_PDE_MOL}) can be 
approximated using finite differences applied to the solution $u$ at the 
current time.
The last term can be computed using either a finite difference approximation 
or a direct analytic calculation.  With the addition of these correction 
terms, the global error can be determined using the analysis in 
Section~\ref{sec:error_analysis}.

Some care must be exercised when choosing the finite difference 
stencils for the correction terms.  While the order of accuracy will always
be boosted, the finite difference stencils used for the correction terms
can introduce errors into the numerical solution that will cause the order
of accuracy to fall below the estimate given by 
(\ref{eq:global_error_ots_simplified}).

\subsubsection*{Semilinear PDEs}
OTS-NIDC can also be applied to semilinear PDEs of the form
\beq
  \frac{\partial u}{\partial t} = A \frac{\partial^n u}{\partial x^n}
  + F \left( \frac{\partial^{n-1} u}{\partial x^{n-1}},
      \frac{\partial^{n-2} u}{\partial x^{n-2}}, \ldots,
      \frac{\partial u}{\partial x}, u \right)
  + f(x,t)
  \label{eq:semilinear_PDE}
\eeq
where $F$ is an arbitrary function of the lower-order spatial derivatives
of $u$ and $f(x,t)$ is a source term.  While the coefficient on the 
leading-order spatial derivative for a semilinear PDE is generally allowed 
to be a function of space and time, optimal time step selection can only be 
applied when it is a constant.  Because the PDE is linear in the leading-order 
spatial derivative, it is natural to construct a finite difference scheme that 
treats $F$ in an explicit manner and lets stability considerations determine 
whether the $A \frac{\partial^n u}{\partial x^n}$ term is handled explicitly or 
implicitly.  Assuming that sufficiently high-order stencils are used to 
compute $F$, the leading-order discretization error for such a finite 
difference scheme is given by
\beq
  \dt \left[ (L u) (\alpha \dx^r - \beta \dt^s) 
  + \left( H u \right) \dx^r 
  + \left( H' u \right) \dt^s
  + (G f) \dt^s \right]
  \label{eq:leading_order_error_semilinear_PDE_MOL}
\eeq
where $L$, $u$, $\alpha$, $\beta,$ $r$, and $s$ are defined as in
(\ref{eq:leading_order_error_model_PDE_general}) and
(\ref{eq:leading_order_error_model_PDE_MOL}), $H$ and $H'$ are spatial 
differential operators (potentially nonlinear) associated with the 
leading-order terms in the discretization error that arise from the nonlinear 
term in~(\ref{eq:semilinear_PDE}), and $G$ is a spatio-temporal differential 
operator that acts on $f$.  With the discretization error written in the 
form (\ref{eq:leading_order_error_semilinear_PDE_MOL}), the optimal time 
step and defect correction terms can be derived by proceeding in the same 
manner as in the previous section.

\subsection{Boundary Conditions}
In general, high-order discretizations of boundary conditions are required 
to ensure that errors introduced at the domain boundary do not pollute 
the numerical solution in the domain interior.  For Dirichlet boundaries, 
it is sufficient to impose boundary conditions by using boundary values 
that have the same order of accuracy as the solution.  More effort is 
required to ensure that Neumann boundary conditions do not reduce the order 
of accuracy of the solution.  In general, Neumann boundary conditions require 
that the accuracy of ghost cells slightly exceeds the desired order of 
accuracy for the solution.  A detailed discussion of methods for imposing 
boundary conditions is beyond the scope of this article.  We refer interested 
readers to the thorough analysis presented in~\cite{gko_book}.  To avoid 
complicating the discussion in Sections~\ref{sec:applications_1d} 
and~\ref{sec:applications_multidim}, Neumann boundary conditions for
all examples are imposed by taking advantage of knowledge of the analytical 
solution.

\subsection{\label{sec:ots_nidc_higher_spatial_dims}
            OTS-NIDC Selection in Multiple Space Dimensions}
OTS-NIDC is applicable to PDEs in any number of space dimensions.  While 
problems in higher dimensions generally require a more careful choice of the 
finite difference stencil used to approximate spatial derivatives, there are 
no fundamental barriers that prevent the use of OTS-NIDC.  Moreover, by using 
a `ghost cell' or `ghost value' method for imposing boundary 
conditions~\cite{gibou_2005,ito_2005,fedkiw_1999,osher_fedkiw_book},
optimal time step selection can even be directly applied to problems on 
irregular domains.  

One important issue that arises in multiple space dimensions is the choice
of grid spacing.  In general, the grid spacing need not be equal in all 
spatial directions.  For some problems, the ratio of the grid spacing in
different coordinate directions may need to be carefully chosen for an 
optimal time step to exist.  In other words, it may be necessary to optimally 
choose both the time step \emph{and} the mesh in order to boost the order 
of accuracy.

\subsection{Limitations of Optimal Time Step Selection}
As with any numerical method, optimal time step selection has its limitations.
The most serious limitation is the requirement that the coefficient on the 
leading-order spatial derivative in the PDE be constant in time and space.  
Without this assumption, the first term in the leading-order discretization 
error would be of the form $(L u) (\alpha(x,t) \dx^r - \beta(x,t) \dt^s) \dt$,
where $\alpha$ and $\beta$ are now functions of space and time.  This
small change in the form of the error makes it impossible to derive a single 
optimal time step.  Instead, each grid point would possess its own individual 
optimal time step.  For a similar reason, optimal time step selection also 
fails for quasilinear and fully nonlinear PDEs.  The issue for these types of 
PDEs is that the coefficients in $P(\dx, \dt)$ become functions of the 
solution $u$.

Another limitation of optimal time step selection is that it does not,
in its current incarnation, apply to finite difference schemes for general
systems of PDEs.  The main issue here is that each PDE might have its 
own optimal time step.  OTS is, however, certainly applicable for problems 
where the optimal time step is the same for all of the individual PDEs in the 
system.

While these limitations certainly place restrictions on the range of problems 
that we can apply OTS-NIDC to, there are many relevant and important PDEs 
arising in practice that satisfy the requirements for OTS-NIDC.

\section{\label{sec:applications_1d} 
         Application to PDEs in One Space Dimension}
Because there are no geometry or anisotropy considerations, OTS-NIDC is very 
useful for boosting the accuracy of many finite difference schemes in one 
spatial dimension.  In this section, we demonstrate its utility in designing 
finite difference schemes that solve several classical linear and semilinear 
PDEs.  To emphasize the importance of including numerical parameter 
optimization as part of the design of any finite difference scheme, we 
include applications of OTS-NIDC to several different types of finite 
difference schemes.  The examples in this section were selected because they 
have analytical solutions that are useful to compare numerical solutions 
against.  However, the utility of OTS-NIDC is by no means limited to 
these simple PDEs.

\subsection{Unit CFL Conditions for Wave Equations}
In this section, we show how OTS selection naturally leads to unit CFL 
conditions for the linear advection equation and second-order wave equation 
in one space dimension.  

\subsubsection{Linear Advection Equation}
When the linear advection equation, $u_t + A u_x = 0$,
is solved using a first-order upwind scheme with forward Euler time
integration, it is well-known that choosing $\Delta t = \Delta x / |A|$ 
completely eliminates all numerical error in the solution (except for 
errors in the initial conditions).  This result is known as the unit CFL 
condition~\cite{leveque_book_2002} and is usually proved by examining how 
characteristic lines of the advection equation pass through grid points.  
Interestingly, the same optimal choice of time step arises by applying OTS 
selection.  

Using OTS selection, we begin by analyzing the discretization error for the 
first-order upwind, forward Euler scheme:
\beq
  u^{n+1}_j = u^{n}_j 
  - A \dt \left( \frac{u^{n}_j - u^{n}_{j-1}}{\dx} \right),
  \label{eq:advection_eqn_1d_FD_scheme}
\eeq
where we have taken $A$ to be positive for convenience.  In general, we only 
need to keep track of the leading-order terms in the discretization error.  
However, for the advection equation, it is valuable to retain all of them.  
Using Taylor series expansions of the true solution of the PDE, we find that 
\bea
  \tu^{n+1}_j &=& \tu^{n}_j 
  + \sum_{k=1}^\infty \frac{\dt^k}{k!} 
       \frac{\partial^k \tu^n}{\pt^k} 
  = \tu^{n}_j + \sum_{k=1}^\infty \frac{\left( -A \dt \right)^k}{k!} 
       \frac{\partial^k \tu^n}{\px^k},
  \label{eq:advection_eqn_1d_time_err} 
  \\
  \frac{\tu^{n}_j - \tu^{n}_{j-1}}{\dx} &=& 
  \frac{\partial \tu^n}{\px} 
  - \frac{1}{\dx} \sum_{k=2}^\infty \frac{\left( -\dx \right)^k}{k!} 
       \frac{\partial^k \tu^n}{\px^k}.
  \label{eq:advection_eqn_1d_space_err}
\eea
where $\tu$ denotes the true solution of the PDE, all derivatives are 
evaluated at $x_j$, and we have used the PDE to replace all of the time 
derivatives with spatial derivatives in (\ref{eq:advection_eqn_1d_time_err}).
For convenience, we will continue to suppress the subscripts on spatial
derivatives throughout the article; the location where they should be 
evaluated will be apparent from the context.

Combining~(\ref{eq:advection_eqn_1d_time_err}) 
and~(\ref{eq:advection_eqn_1d_space_err}) with the finite difference 
scheme~(\ref{eq:advection_eqn_1d_FD_scheme}) yields the evolution
equation for the error $e \equiv u - \tu$:
\beq
  e^{n+1}_j = e^{n}_j 
    - A \dt \left( \frac{e^{n}_j - e^{n}_{j-1}}{\dx} \right) 
    + A \dt \sum_{k=2}^\infty \frac{\left( -1 \right)^k}{k!} 
        \frac{\partial^k \tu^n}{\px^k} 
        \left( \dx^{k-1} - \left( A \dt \right)^{k-1} \right).
  \label{eq:advection_eqn_1d_err_eqn}
\eeq
Therefore, choosing the time step to eliminate the leading-order term 
of the discretization error leads to the unit CFL condition and 
\emph{complete elimination} of the local truncation error at each time step.
It is interesting to note that OTS selection applied to the Lax-Wendroff 
scheme~\cite{leveque_book_2002,leveque_book_1992} yields the same optimal 
time step and total elimination of the truncation error.

\subsection{Second-Order Wave Equation\label{sec:wave_eqn_1d}}
Numerical methods for the wave equation $u_{tt} - c^2 u_{xx} = f$ are typically 
derived by transforming the second-order equation into an equivalent
first-order system of equations.  However, as Kreiss, Petersson, and Ystr\"om 
showed in~\cite{kreiss2002}, direct discretization of the second-order wave 
equation is a viable option and leads to the following finite difference 
scheme that is second-order accurate in both space and time:
\bea
  \frac{u^{n+1}_i - 2 u^n_i + u^{n-1}_i}{\dt^2}
  = c^2 \left( \frac{u^{n}_{i+1} - 2 u^n_i + u^n_{i-1}}{\dx^2} \right)
  + f.
  \label{eq:wave_eqn_KPY}
\eea
The stability constraint for this scheme of the form $\dt = K \dx$.

\begin{figure}[tb]
\begin{center}
\scalebox{0.35}{\includegraphics{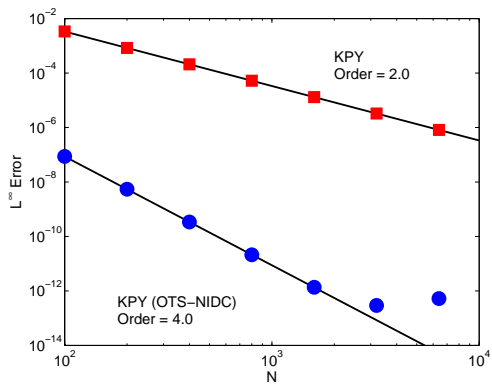}}
\caption{$L^\infty$ error as a function of number of grid points for the
  KPY discretization of the second-order wave equation on with OTS-NIDC 
  (blue circles) and without OTS-NIDC (red squares).
}
\label{fig:wave_eqn_1d_error}
\end{center}
\end{figure}

To compute the optimal time step and correction terms for this scheme, we
begin by deriving the truncation error for the scheme.  Employing Taylor 
series expansions, it is straightforward to show that the true solution 
satisfies
\bea
  \frac{\tu^{n+1}_i - 2 \tu^n_i + \tu^{n-1}_i}{\dt^2}
    - \frac{\dt^2}{12} \frac{\partial^4 \tu}{\pt^4}
  &=& c^2 \left( \frac{\tu^{n}_{i+1} - 2 \tu^n_i + \tu^n_{i-1}}{\dx^2}
  -\frac{\dx^2}{12} \frac{\partial^4 \tu}{\px^4} \right) + f
  \nonumber \\
  &+& O(\dt^4) + O(\dx^4) 
\eea
Combining this result with the observation that
$\tu_{tttt} = c^4 \tu_{xxxx} + c^2 f_{xx} + f_{tt}$ and rearranging a bit, 
we find that the error equation for (\ref{eq:wave_eqn_KPY}) is given by
\bea
  e^{n+1}_i &=& 2 e^n_i - e^{n-1}_i
  + \dt^2 c^2 \left( \frac{e^{n}_{i+1} - 2 e^n_i + e^n_{i-1}}{\dx^2} \right)
  + \frac{\dt^2}{12} c^2 \frac{\partial^4 \tu}{\px^4} 
    \left( \dx^2 - \dt^2 c^2 \right)
  \nonumber \\
  &-& \frac{\dt^4}{12} \left( c^2 \frac{\partial^2 f}{\px^2} 
                            + \frac{\partial^2 f}{\pt^2} \right)
      + O(\dt^2 \dx^4) + O(\dt^6).
\eea
Choosing $\dt = \dx/c$ and adding the correction term 
$\dt^2 \left(c^2 f_{xx} + f_{tt} \right)/12$ to the right-hand side of 
(\ref{eq:wave_eqn_KPY}) eliminates the leading-order term in the local 
truncation error.  Using the heuristic for two-step methods that the global 
error should be approximately $1/\dt^2$ times the local error leads to a 
global error of $O(\dx^4) = O(\dt^4)$ -- a boost in the order of accuracy 
from second- to fourth-order.  Notice that as for first-order upwind, forward 
Euler and Lax-Wendroff schemes for the linear advection equation, applying OTS 
selection to the KPY scheme leads to an optimal time step that satisfies the 
unit CFL condition.  However, unlike the schemes for the linear advection 
equation, the unit CFL condition for the KPY scheme eliminates only the 
leading-order term in the truncation error.  Figure~\ref{fig:wave_eqn_1d_error} 
demonstrates that the expected accuracy is indeed achieved by applying 
OTS-NIDC to the KPY scheme.  

An important property of the KPY scheme is that the error introduced during 
the first time step affects the global error at all times.  More specifically, 
a $(p+1)$-order accurate first time introduces an $O(\dx^p)$ error into the 
solution at all future times.
Thus, we must ensure that the error introduced by the first time step is
at least one order of accuracy higher than desired accuracy for the numerical 
solution.  See Appendix~\ref{appendix:KPY_analysis} for a derivation of this
condition. 

Constructing a fifth-order approximation for the first time step is 
straightforward using a fifth-order Taylor series expansion in time:
\bea
  u^1 &=& u^0 + \dt \frac{\partial u^0}{\pt} 
              + \frac{\dt^2}{2} \frac{\partial^2 u^0}{\pt^2}
              + \frac{\dt^3}{6} \frac{\partial^3 u^0}{\pt^3} 
           + \frac{\dt^4}{24} \frac{\partial^4 u^0}{\pt^4}
  \nonumber \\
  &=& u^0 + \dt \frac{\partial u^0}{\pt} 
          + \frac{\dt^2}{2} \left(c^2 \frac{\partial^2 u^0}{\px^2} + f\right)
  + \frac{\dt^3}{6} \left(c^2 
    \frac{\partial^2}{\px^2} \left(\frac{\partial u^0}{\pt} \right) 
    + \frac{\partial f}{\pt}\right)
  \nonumber \\
  & & +\  \frac{\dt^4}{24} \left(c^4 \frac{\partial^4 u^0}{\px^4} 
  + c^2 \frac{\partial ^2 f}{\px^2} + \frac{\partial^2 f}{\pt^2} \right).
  \label{eq:KPY_fifth_order_first_step}
\eea
Note that special care must be exercised if finite differences are used to 
compute the derivatives in (\ref{eq:KPY_fifth_order_first_step}) because
higher-order terms may be implicitly included by lower-order terms in the 
expansion if the initial time step is taken using $\dto$.

\subsection{Diffusion Equation \label{sec:diffusion_eqn_1d}}
Solution of the diffusion equation, $u_t = D u_{xx}$, using forward Euler 
time integration with a second-order central difference approximation for 
the Laplacian provides a familiar example where careful choice of the time 
step leads to cancellation of the leading order spatial and temporal 
truncation errors.  It is well-known that choosing the time step to be 
$\dto = \dx^2/6D$ boosts the order of accuracy for this simple scheme 
from $O(\dx^2)$ to $O(\dx^4)$.

When a source term is present, boosting the order of accuracy requires
the addition of defect correction terms.  These correction terms are 
straightforward to calculate by examining the discretization error for the
scheme
\beq
  u^{n+1}_j = u^{n}_j 
  + \dt 
    \left( D \left[\frac{u^{n}_{j+1} -2 u^{n}_j + u^{n}_{j-1}}{\dx^2}\right] 
         + f^n_j \right),
  \label{eq:diffusion_eqn_1d_FD_scheme}
\eeq
where $f$ is the source term.  Note that even though this scheme is formally 
first-order in time and second-order in space, the stability constraint 
$\dt \le \dx^2/2D$ implies that the scheme is $O(\dx^2)$ accurate 
overall (\ie inclusion of the temporal order of accuracy is irrelevant).
Using Taylor series expansions and the PDE $u_t = D u_{xx} + f$,
we see that the true solution satisfies
\bea
  \tu^{n+1}_j = \tu^{n}_j 
  &+& \dt \left( D \frac{\partial^2 \tu^n}{\px^2} + f^n_j \right)
  \nonumber \\
  &+& \frac{(D\dt)^2}{2} \frac{\partial^4 \tu^n}{\px^4} 
  + \frac{\dt^2}{2} \left( D\frac{\partial^2 f}{\px^2}
                         + \frac{\partial f}{\pt} \right)
  + O \left( \dt^3 \right)
  \label{eq:diffusion_eqn_1d_time_err} 
\eea
and that the central difference approximation for the Laplacian satisfies
\beq
  \frac{\tu^{n}_{j+1} -2 \tu^{n}_j + \tu^{n}_{j-1}}{\dx^2}  =
  \frac{\partial^2 \tu^n}{\px^2} 
  + \frac{\dx^2}{12} \frac{\partial^4 \tu^n}{\px^4} 
  + O(\dx^4).
  \label{eq:diffusion_eqn_1d_space_err}
\eeq
Therefore, the truncation error for (\ref{eq:diffusion_eqn_1d_FD_scheme})
is given by
\beq
  \frac{\partial^4 \tu^n}{\px^4} 
    \left[ \frac{\dx^2}{12} - \frac{D \dt}{2} \right] (D \dt)
    - \frac{\dt^2}{2} \left( D \frac{\partial^2 f}{\px^2} 
                           + \frac{\partial f}{\pt} \right)
      + O(\dt \dx^4) + O(\dt^3).
  \label{eq:diffusion_eqn_1d_trunc_err}
\eeq
From this expression, we easily recover the well-known optimal time step for
the diffusion equation, $\dto = \dx^2/6D$, and identify the correction term 
$\dt^2 \left( D f_{xx} + f_t \right)  /2$.

\begin{figure}[tb]
\begin{center}
\scalebox{0.35}{\includegraphics{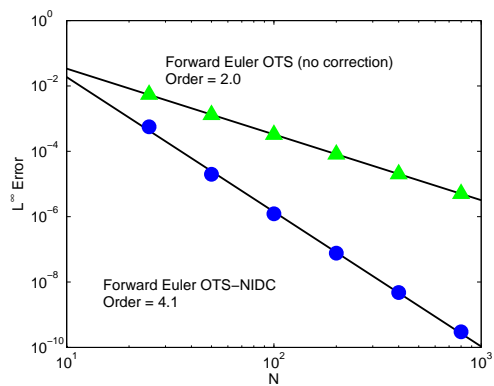}} 
\caption{$L^\infty$ error as a function of number of grid points for 
two finite difference schemes that solve the diffusion equation: 
forward Euler with optimal time step and defect correction term (circles) and 
forward Euler with optimal time step \emph{without} defect correction term 
(triangles).  
}
\label{fig:diffusion_eqn_1d_src_error}
\end{center}
\end{figure}

From (\ref{eq:global_error_ots_simplified}), we see that using the optimal 
time step and adding this simple correction term to 
(\ref{eq:diffusion_eqn_1d_FD_scheme}) boosts the accuracy of original scheme 
from second to fourth order.  Note that to achieve fourth order accuracy, 
the spatial derivatives in the correction term must be calculated using a
finite difference approximation that is at least second-order accurate; 
otherwise, the truncation error will be $O(\dt^2 \dx)$ which implies an 
$O(\dt \dx) = O(\dx^3)$ global error. 
The importance of the correction term is illustrated in 
Figure~\ref{fig:diffusion_eqn_1d_src_error}, which compares the accuracy of 
the OTS forward Euler scheme with and without the correction term.  As we 
can see, OTS forward Euler with correction term is fourth-order accurate 
while OTS forward Euler without correction term is only second-order 
accurate.  

\subsubsection{DuFort-Frankel Scheme \label{sec:dufort_frankel}}
As mentioned earlier, the forward Euler scheme for the diffusion equation is
unstable for problems in more than 3 space dimensions.  Fortunately, the 
OTS-NIDC makes it possible to construct a fourth-order accurate 
DuFort-Frankel scheme. 

We begin by writing the DuFort-Frankel scheme in form that is easily 
generalizable to higher space dimensions:
\beq
  u^{n+1}_j = u^{n-1}_j 
  + 2 \dt 
    \left( D \left[\frac{u^{n}_{j+1} -2 u^{n}_j + u^{n}_{j-1}}{\dx^2}\right] 
         + f^n_j \right)
  - 2 D \frac{\dt^3}{\dx^2}
    \left( \frac{u^{n+1}_j -2 u^{n}_j + u^{n-1}_j}{\dt^2}\right).
  \label{eq:diffusion_eqn_1d_DF_scheme}
\eeq
Note that the DuFort-Frankel scheme is essentially a leap-frog time integration 
scheme that is stabilized by the addition of a term proportional to $u_{tt}$
that vanishes in the limit $\dt \rightarrow 0$.  This scheme is 
unconditionally stable and has a global error that is formally 
$O(\dt^2/\dx^2) + O(\dt^2) + O(\dx^2) + O(\dt^4/\dx^2)$~\cite{gko_book}.  
This expression for the formal error is ``minimized'' when $\dt = O(\dx^2)$, 
which yields an $O(\dx^2)$ global error.

\begin{figure}[tb]
\begin{center}
\scalebox{0.35}{\includegraphics{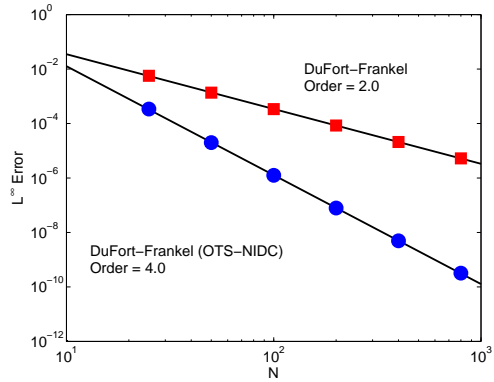}} 
\caption{$L^\infty$ error in the numerical solution of the 1D diffusion 
equation as a function of number of grid points for the DuFort-Frankel scheme 
with (circles) and without (squares) OTS-NIDC.
}
\label{fig:diffusion_eqn_1d_dufort_frankel_error} 
\end{center}
\end{figure}

Using Taylor series expansions centered at $(x_j, t_n)$, 
(\ref{eq:diffusion_eqn_1d_space_err}) and the PDE to derive the local 
truncation error for the DuFort-Frankel scheme, we find that 
\bea
  \tau &=& 2D \dt \frac{\partial^4 \tu}{\px^4}
  \left(\frac{\dx^2}{12} - \frac{D^2 \dt^2}{\dx^2} \right)
  - 2 D \frac{\dt^3}{\dx^2} \left(D \frac{\partial^2 f}{\px^2}  
                                 + \frac{\partial f}{\pt} \right)
  \nonumber \\
  &+& O(\dt^3) + O(\dt \dx^4) + O\left(\frac{\dt^5}{\dx^2} \right).
\eea
Therefore the optimal time step is $\dto = \dx^2/D\sqrt{12}$ and the defect
correction term is $2 D \dt^3 \left(D f_{xx} + f_{t} \right)/\dx^2$.  As with 
the forward Euler scheme for the diffusion equation, OTS-NIDC boosts the
order of accuracy for the DuFort-Frankel scheme from two to four.  
Figure~\ref{fig:diffusion_eqn_1d_dufort_frankel_error} shows that OTS-NIDC
applied to the DuFort-Frankel scheme achieves the expected boost in accuracy. 

\subsection{Viscous Burgers Equation}
We now demonstrate the use of OTS-NIDC in the context of a well-studied
semilinear PDE -- the viscous Burgers equation~\cite{whitham_book}:
\beq
  \frac{\partial u}{\pt} + u \frac{\partial u}{\px} = 
     \nu \frac{\partial^2 u}{\px^2} 
  \label{eq:burgers_1d}
\eeq
where $\nu$ is the viscosity.  For this equation, let us consider the finite
difference scheme constructed using a forward Euler time integration scheme
and second-order, central difference approximations for the diffusion and 
nonlinear advection terms:
\beq
  u^{n+1}_j = u^{n}_j 
  + \dt 
    \left( \nu 
      \left [ \frac{u^{n}_{j+1} -2 u^{n}_j + u^{n}_{j-1}}{\dx^2} \right]
         - u^n_j 
      \left[ \frac{u^{n}_{j+1} - u^{n}_{j-1}}{2 \dx} \right] 
    \right).
  \label{eq:burgers_1d_FD_scheme}
\eeq
This scheme is formally first-order in time and second-order in space.  
Using the stability condition for the advection-diffusion as a heuristic 
guide\footnote{The stability constraint for the advection-diffusion equation
is $\dt < \min \{2D/A^2, \dx^2/2D \}$~\cite{chan_1984}.}, 
we expect a stability constraint of the form $\dt \le \dx^2/2 \nu$ in the 
limit $\dx \rightarrow 0$ (assuming that the solution $u$ remains bounded).  
Therefore, the scheme (\ref{eq:burgers_1d_FD_scheme}) is $O(\dx^2)$ accurate 
overall. 

\begin{figure}[tb]
\begin{center}
\scalebox{0.35}{\includegraphics{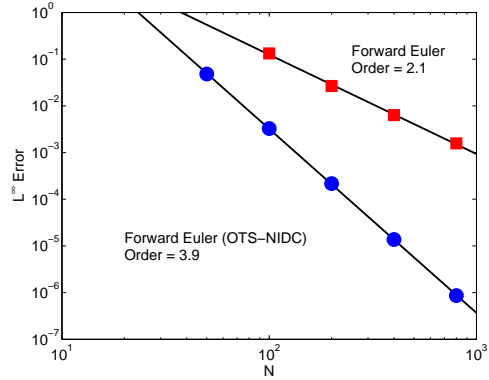}} 
\caption{$L^\infty$ error as a function of number of grid points for
two finite difference schemes that solve the viscous Burgers equation:
forward Euler with OTS-NIDC (circles)
and forward Euler with suboptimal time step $\dt = \dx^2/4\nu$ and no 
correction term (squares).
The viscous Burgers equation is solved with initial conditions taken to be 
$u(x,0) = 1 + \gamma \sqrt{\frac{\nu}{\pi}} \exp\left( -(x-1)^2/4 \nu \right) 
            / \left[ 1 + \frac{\gamma}{2} \ 
                         \mathtt{erfc}\left((x-1)/\sqrt{4\nu} \right) 
    \right]$,
where $\ln(1+\gamma)$ is the effective Reynolds number (set equal to 10) and 
the wave speed is set equal to $1$.  The boundary conditions are imposed 
exactly by using
the analytical solution
$u(x,t) = 1 + \gamma \sqrt{\frac{\nu}{\pi T}} 
              \exp\left( -(x-T)^2/4 \nu T \right) 
            / \left[ 1 + \frac{\gamma}{2} \ 
                         \mathtt{erfc}\left((x-T)/\sqrt{4\nu T} \right) 
    \right]$,
where $T = t+1$.
}
\label{fig:burgers_1d_error}
\end{center}
\end{figure}

To analyze the discretization error for (\ref{eq:burgers_1d_FD_scheme}), 
we combine~(\ref{eq:diffusion_eqn_1d_space_err}) from our analysis of 
the diffusion equation with the leading-order discretization error for 
the central difference approximation of $u_x$,
\beq
  \frac{\tu^{n}_{j+1} - \tu^{n}_{j-1}}{2 \dx}  =
  \frac{\partial \tu^n}{\px} 
  + \frac{\dx^2}{6} \frac{\partial^3 \tu^n}{\px^3} 
  + O(\dx^4),
  \label{eq:burgers_1d_ux_err}
\eeq
to obtain the local truncation error for~(\ref{eq:burgers_1d_FD_scheme}): 
\bea
  & &
      \left( \nu \frac{\partial^4 \tu^n}{\px^4} 
           - 2 \tu^n_j \frac{\partial^3 \tu^n}{\px^3} \right)
      \left[ \frac{\dx^2}{12} - \frac{\nu \dt}{2}  \right] \dt
  \nonumber \\
  &+& \frac{\dt^2}{2} 
      \left( 
           4 \nu \frac{\partial \tu^n}{\px} \frac{\partial^2 \tu^n}{\px^2}
         - 2 \tu^n_j \left( \frac{\partial^2 \tu^n}{\px^2} \right)^2
         - \left(\tu^n_j\right)^2 \frac{\partial^2 \tu^n}{\px^2}
      \right) 
  \nonumber \\
  &+& O(\dt^3) + O(\dt \dx^4).
  \label{eq:burgers_1d_err_eqn}
\eea
Therefore, the optimal time step is $\dto = \dx^2/6\nu$ and the defect 
correction term is
\beq
  - \frac{\dt^2}{2} 
      \left( 
           4 \nu \frac{\partial \tu^n}{\px} \frac{\partial^2 \tu^n}{\px^2}
         - 2 \tu^n_j \left( \frac{\partial^2 \tu^n}{\px^2} \right)^2
         - \left(\tu^n_j\right)^2 \frac{\partial^2 \tu^n}{\px^2}
      \right)
  \label{eq:burgers_1d_corr_term}
\eeq 
Note that the optimal time step eliminates errors introduced by 
both the diffusion term \emph{and} the nonlinear advection term.

By using the optimal time step and the correction term, the overall error
for the finite difference scheme~(\ref{eq:burgers_1d_FD_scheme}) is reduced
from $O(\dx^2)$ to $O(\dx^4)$.  As for the diffusion equation, we need
to use a second-order accurate discretization for the spatial derivatives
in the correction term.  
Figure~\ref{fig:burgers_1d_error} demonstrates the boost in accuracy obtained 
by applying OTS-NIDC to the finite difference 
scheme~(\ref{eq:burgers_1d_FD_scheme}).

The Burgers equation example illustrates an important aspect of OTS-NIDC:
the ease with which OTS-NIDC can be used to boost the order of accuracy 
depends on the choice of the finite difference scheme.  Had we chosen an 
upwind approximation for the gradient term in the Burgers equation, the 
correction term would have been more complicated and included terms involving 
higher-order spatial derivatives than those that appear in the original 
equation.  While it should be possible to handle this situation with OTS-NIDC 
selection, it is clearly less desirable.

\subsection{Fourth-Order Parabolic Equation}
As a final example in one space dimension, we consider a PDE with high-order
spatial derivatives -- a fourth-order parabolic equation: 
\beq
  \frac{\partial u}{\pt} = -\kappa \frac{\partial^4 u}{\px^4} + f(x,t), 
  \label{eq:4th_order_parabolic_eqn_1d}
\eeq
where $\kappa > 0$.  This example illustrates the importance of the choice of
discretization for the spatial derivative in order for OTS-NIDC to be
applicable. 

From our previous examples, we know that optimal time step selection requires 
the coefficient on leading-order spatial discretization error to be directly 
related to the second derivative of the solution $u$ with respect to time:
\beq
  \frac{\partial^2 u}{\pt^2} = 
    \kappa^2 \frac{\partial^8 u}{\px^8} 
  - \kappa \frac{\partial^4 f}{\px^4} 
  + \frac{\partial f}{\pt}
  \label{eq:4th_order_parabolic_eqn_1d_second_time_derivative}.
\eeq
Therefore, the leading-order term in the spatial discretization 
error should involve the eighth-derivative of $u$, which implies that the
finite difference approximation to the bilaplacian needs to be at least 
fourth-order accurate.  With this insight, let us consider a finite 
difference scheme for~(\ref{eq:4th_order_parabolic_eqn_1d}) that 
uses forward Euler for time integration and a fourth-order central 
difference approximation for the bilaplacian:
\beq
  u^{n+1}_j = u^{n}_j + \dt \left( -\kappa B u^n + f_j^n \right), 
  \label{eq:4th_order_parabolic_eqn_1d_FD_scheme}
\eeq
where the fourth-order discrete bilaplacian is given by
\beq
  B u^n = \frac{ -u^{n}_{j+3} + 12 u^{n}_{j+2} - 39 u^{n}_{j+1}
               + 56 u^{n}_j
               - 39 u^{n}_{j-1} + 12 u^{n}_{j-2} -u^{n}_{j-3} }
               {6 \dx^4}.
\eeq
This scheme is formally first-order in time and fourth-order in space.  
Due to the stability constraint $\dt \le 3\dx^4/40 \kappa$, the scheme is 
$O(\dx^4)$ accurate overall.

\begin{figure}[tb]
\begin{center}
\scalebox{0.33}{
  \includegraphics{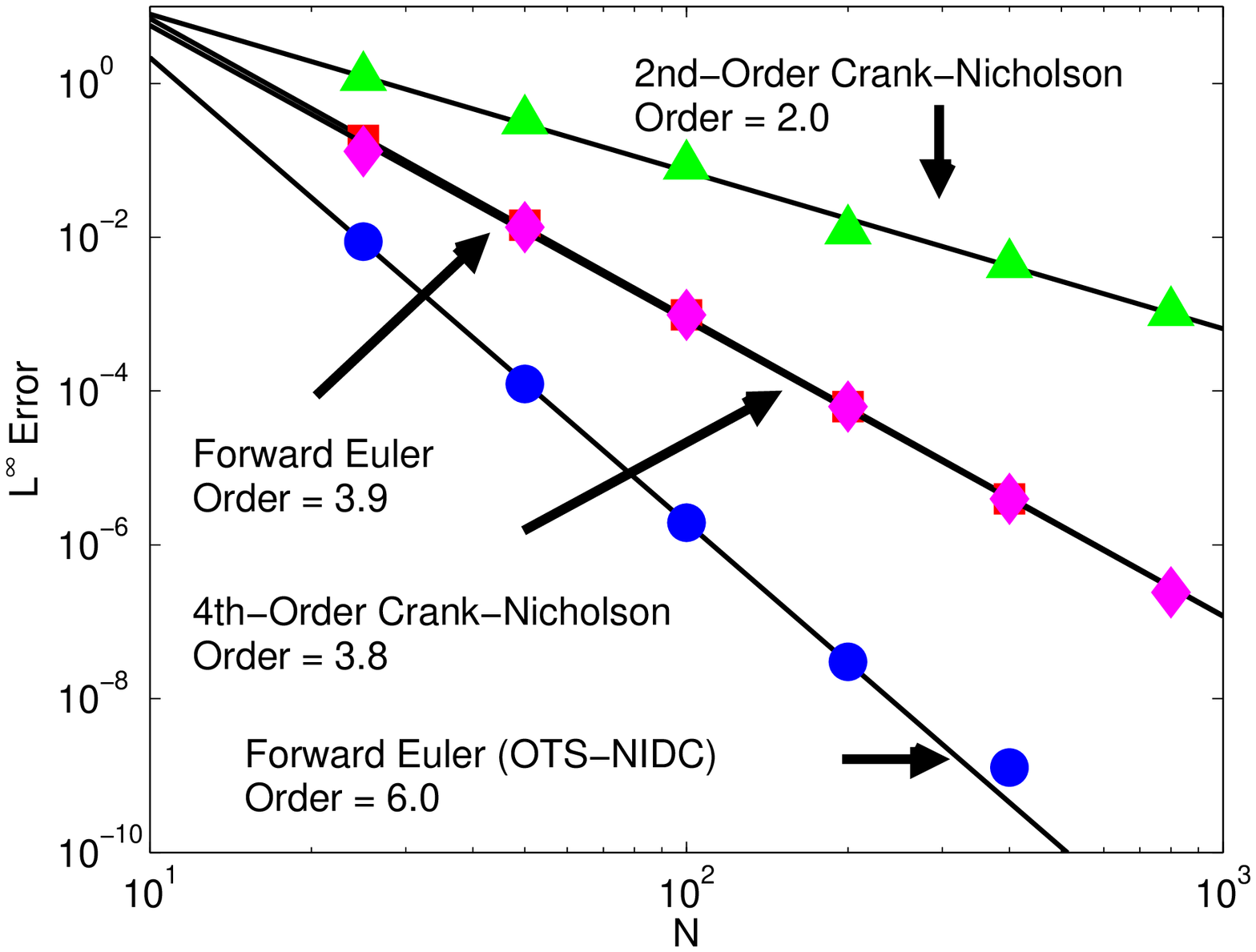}} 
\ \ 
\scalebox{0.33}{
  \includegraphics{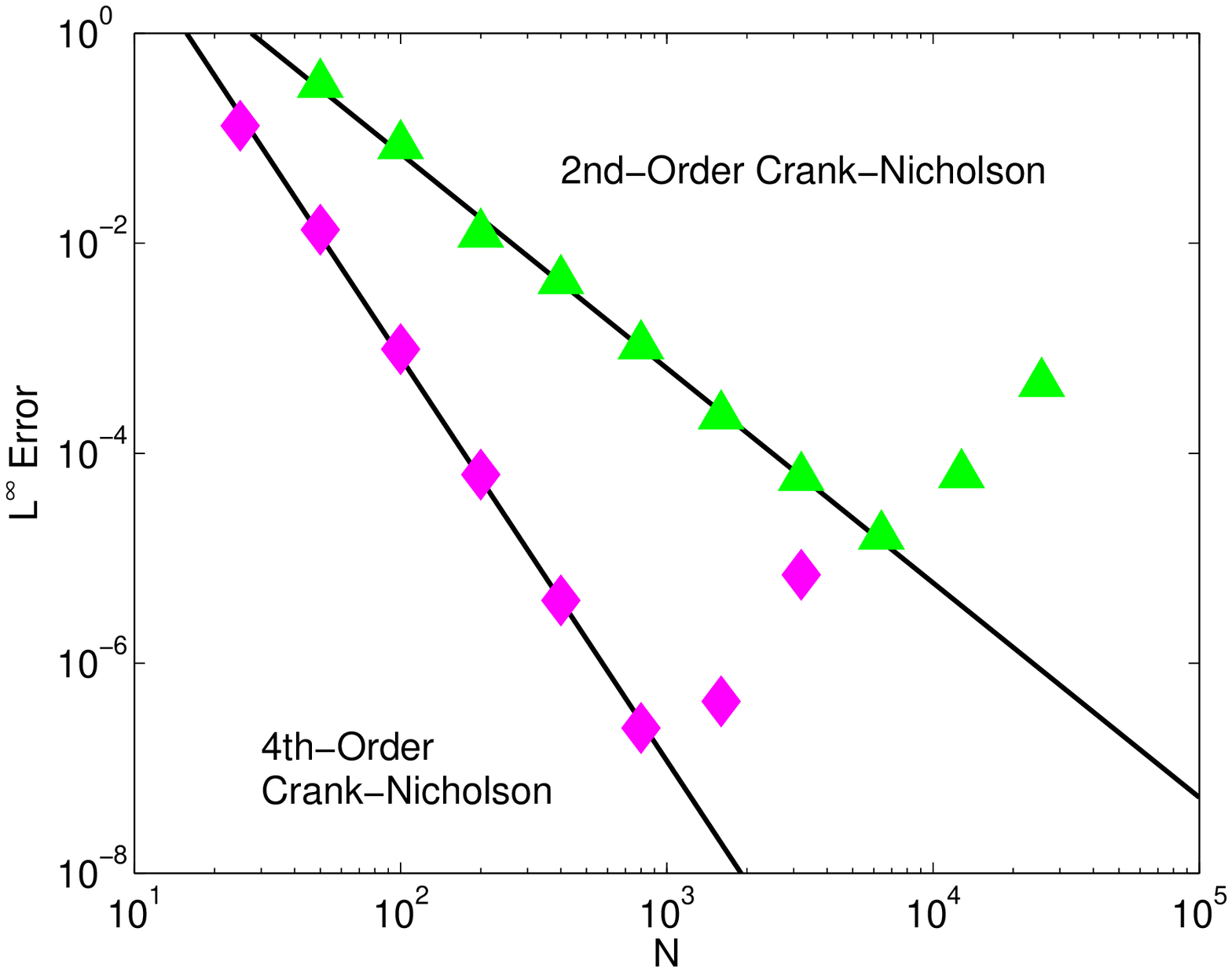}} 
\caption{$L^\infty$ error as a function of number of grid points for four
finite difference schemes that solve the fourth-order parabolic equation
(\ref{eq:4th_order_parabolic_eqn_1d}):
forward Euler with OTS-NIDC (circles),
forward Euler with suboptimal time step and no correction term (squares),
Crank-Nicholson with a fourth-order stencil for the bilaplacian and
an $O(\dx^2)$ time step (diamonds), and Crank-Nicholson with a second-order 
stencil for the bilaplacian (triangles).
Both forward Euler schemes use a fourth-order, central finite difference 
approximation for the bilaplacian.  The boundary and ghost cells values
are set by taking values from the known analytical solution.
Note that the errors for the forward Euler with suboptimal time step and 
fourth-order Crank-Nicholson schemes lie almost directly on top of each 
other at the resolution of the figures.  
The figure on the left highlights the impact (when $N$ is large) of round-off 
error associated with matrix inversion on the accuracy of the solution for the 
two Crank-Nicholson schemes.
}
\label{fig:4th_order_parabolic_eqn_1d_error}
\end{center}
\end{figure}

Proceeding in the usual way, we find that the true solution satisfies
\beq
  \tu^{n+1}_j = \tu^{n}_j 
  + \dt \left( -\kappa \frac{\partial^4 \tu^n}{\px^4}  
             + f^n_j
        \right)
  + \frac{\dt^2}{2} 
    \left(
      \kappa^2 \frac{\partial^8 \tu^n}{\px^8} 
    - \kappa \frac{\partial^4 f}{\px^4} 
    + \frac{\partial f}{\pt}
    \right)
  + O \left( \dt^3 \right)
  \label{eq:4th_order_parabolic_eqn_1d_time_err}
\eeq
and that the central difference approximation to the bilaplacian satisfies
\beq
  B \tu^n = \frac{\partial^4 \tu}{\px^4} 
  - \frac{7 \dx^4}{240} 
    \frac{\partial^8 \tu}{\px^8} 
  + O(\dx^6)
  \label{eq:4th_order_parabolic_eqn_1d_space_err}.
\eeq
Combining these results with the finite difference 
scheme~(\ref{eq:4th_order_parabolic_eqn_1d_FD_scheme}), we find that the
local truncation error is given by
\beq
     \frac{\partial^8 \tu^n}{\px^8} 
    \left[ \frac{7 \dx^4}{240} - \frac{\kappa \dt}{2}  \right] (\kappa \dt)
    + \frac{\dt^2}{2} 
        \left( \kappa \frac{\partial^4 f}{\px^4} 
             - \frac{\partial f}{\pt}
        \right)
    + O(\dt^3) + O(\dt \dx^6)
  \label{eq:4th_order_parabolic_eqn_1d_err_eqn}
\eeq
Therefore, the optimal time step is $\dto = 7 \dx^4/120\kappa$ and the 
correction term is 
\beq
  - \frac{\dt^2}{2} 
      \left( \kappa \frac{\partial^4 f}{\px^4} 
           - \frac{\partial f}{\pt}
      \right).
\eeq
Together, these modifications boost the accuracy of the original
finite difference scheme from fourth- to sixth-order.
As with the previous examples, the spatial derivatives in the 
correction term should be computed using a finite difference approximation 
that is at least second-order accurate.  
The improved accuracy of (\ref{eq:4th_order_parabolic_eqn_1d_FD_scheme}) 
when the optimal time step and defect correction terms are used is 
illustrated in Figure~\ref{fig:4th_order_parabolic_eqn_1d_error}. 

Applying optimal time step selection to high-order PDEs requires the use of 
high-order discretizations of the spatial derivatives.  While the need to 
derive high-order finite difference stencils may seem like an unnecessary 
burden, it is important to remember that high-order finite difference schemes 
for spatial derivatives are almost always necessary to maximize the efficiency 
of explicit time integration schemes for high-order 
PDEs\footnote{Because restrictive time step constraints are always a major 
limitation of explicit schemes for time dependent PDEs with high-order spatial 
derivatives~\cite{gko_book,greer_2006}, 
the overall error of finite difference schemes for high-order PDEs is usually 
dominated by the error introduced by the spatial discretization.  As a result, 
it is crucial that numerical schemes for high-order PDEs use high-order 
discretizations for the spatial derivatives.}.
Optimal time step selection merely provides a valuable way to leverage the 
effort spent deriving high-order finite difference schemes in order to further 
boost the overall accuracy of the numerical method.

\begin{table}[tb]
\caption{
Computational cost as a function of the global numerical error $e$
for various finite difference schemes that numerically solve the 1D 
fourth-order parabolic equation.  For all finite difference schemes, the 
time step is chosen so that the temporal error is subdominant to the spatial
error.
}
\label{tab:comp_perf_vs_err_4th_order_parabolic} 
\renewcommand{\arraystretch}{1.3}
\centering
\begin{tabular}{lcccc}
  \hline
  {\bf Numerical Scheme} & $\dx$ 
  & $\dt$
  & {\bf Memory}
  & {\bf Compute Time}
  \\
  \hline 
  Forward Euler with 
    & \multirow{2}{*}{$O\left( e^{1/2} \right)$} 
    & \multirow{2}{*}{$O\left( e^{2} \right)$}
    & \multirow{2}{*}{$O\left( e^{-1/2} \right)$} 
    & \multirow{2}{*}{$O\left( e^{-5/2} \right)$} \\
  \ \ 2nd-Order Bilaplacian & & & & \\
  Forward Euler with 
    & \multirow{2}{*}{$O\left( e^{1/4} \right)$}
    & \multirow{2}{*}{$O\left( e \right)$}
    & \multirow{2}{*}{$O\left( e^{-1/4} \right)$} 
    & \multirow{2}{*}{$O\left( e^{-{5/4}} \right)$} \\
  \ \ 4th-Order Bilaplacian & & & & \\
  Forward Euler with 
    & \multirow{2}{*}{$O\left( e^{1/6} \right)$} 
    & \multirow{2}{*}{$O\left( e^{2/3} \right)$} 
    & \multirow{2}{*}{$O\left( e^{-1/6} \right)$} 
    & \multirow{2}{*}{$O\left( e^{-5/6} \right)$} \\ 
  \ \ 4th-Order Bilaplacian and OTS-NIDC & & & & \\
  Crank-Nicholson with 
    & \multirow{2}{*}{$O\left( e^{1/2} \right)$} 
    & \multirow{2}{*}{$O\left( e^{1/2} \right)$} 
    & \multirow{2}{*}{$O\left( e^{-1/2} \right)$} 
    & \multirow{2}{*}{$O\left( e^{-1} \right)$} \\
  \ \ 2nd-Order Bilaplacian & & & & \\
  Crank-Nicholson with 
    & \multirow{2}{*}{$O\left( e^{1/4} \right)$} 
    & \multirow{2}{*}{$O\left( e^{1/2} \right)$} 
    & \multirow{2}{*}{$O\left( e^{-1/4} \right)$} 
    & \multirow{2}{*}{$O\left( e^{-3/4} \right)$} \\
  \ \ 4th-Order Bilaplacian & & & & \\
  \hline
\end{tabular}
\end{table}

\begin{figure}[tb]
\begin{center}
\scalebox{0.33}{
  \includegraphics{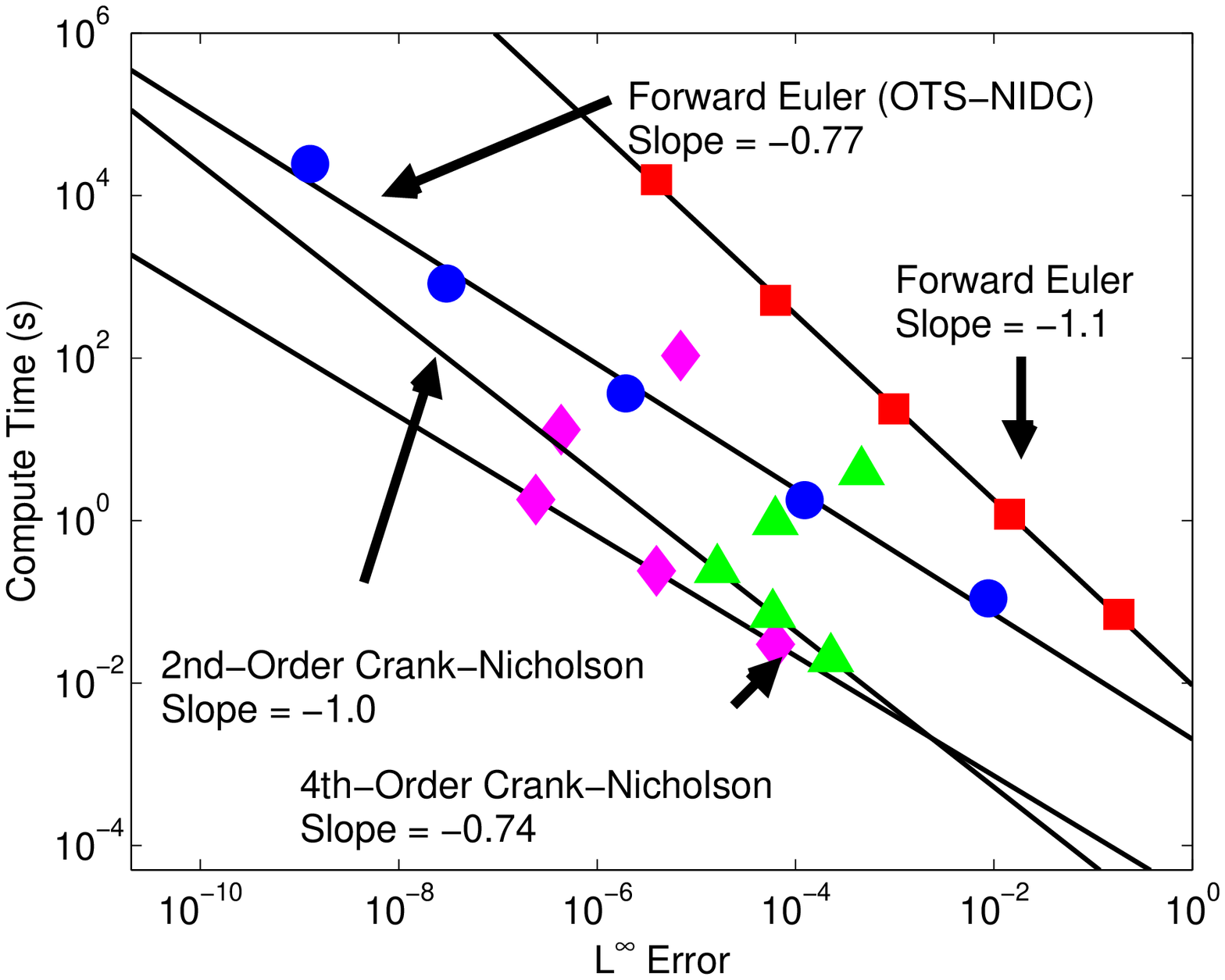}} 
\caption{Computation time as a function of $L^\infty$ error for four
finite difference schemes that solve the fourth-order parabolic equation
(\ref{eq:4th_order_parabolic_eqn_1d}).  
These results were obtained using MATLAB implementations of the 
finite difference schemes run on a 2.4 GHz MacBook Pro.
Note the impact of round-off error associated with matrix inversion on 
the accuracy of the solution for the two Crank-Nicholson schemes.
See the caption of Figure~\ref{fig:4th_order_parabolic_eqn_1d_error} for 
descriptions of the four finite difference schemes.
}
\label{fig:4th_order_parabolic_eqn_1d_perf}
\end{center}
\end{figure}

Table~\ref{tab:comp_perf_vs_err_4th_order_parabolic} shows the computational 
performance of various finite difference methods for solving 
equation~(\ref{eq:4th_order_parabolic_eqn_1d_second_time_derivative}).
Observe the dramatic decrease in computational cost as a function of the
error by simply replacing a second-order discretization for the bilaplacian 
operator with a fourth-order discretization.  For the forward Euler scheme,
using an OTS-NIDC drives the computational cost down even more.
The performance plot in Figure~\ref{fig:4th_order_parabolic_eqn_1d_perf} 
shows that the theoretical estimates are in good agreement with timings 
from computational experiments.  Note that the empirical performance for the 
forward Euler schemes is better than theoretically expected because the 
computational work per time step is not dominated by the 
$O(1/\dx)$ update calculations until $\dx$ gets very small.

In Figure~\ref{fig:4th_order_parabolic_eqn_1d_perf}, observe that the 
2nd-order Crank-Nicholson scheme actually outperforms the OTS-NIDC scheme over 
the high (but still practically useful) range of errors even though 
Table~\ref{tab:comp_perf_vs_err_4th_order_parabolic} predicts that the OTS-NIDC
scheme should do better.  The origin of this discrepancy is the constant 
hidden by the big-$O$ notation.  In particular, the factor of $7/120$ in 
the optimal time step for the forward Euler scheme with a fourth-order 
approximation of the bilaplacian leads to a relatively large hidden constant 
for the compute time. 
However, as Figures~\ref{fig:4th_order_parabolic_eqn_1d_error} and
\ref{fig:4th_order_parabolic_eqn_1d_perf} show, there is an important 
limitation of the 2nd-order Crank-Nicholson scheme even if we ignore the 
cross-over which will occur at sufficiently small values of the error -- 
round-off error associated with matrix inversion limits the accuracy that is 
achievable by the scheme.  This limitation also exists for the 4th-order 
Crank-Nicholson scheme. 

Comparing the computational cost of the best forward Euler and Crank-Nicholson 
schemes, we see that the OTS-NIDC forward Euler scheme more efficiently 
utilizes memory at the cost of requiring more computation time.  This tradeoff 
between memory and computation only exists in one space dimension.  In 
two-dimensions, both schemes require the same amount of compute time; in three 
or more space dimensions, the forward Euler scheme requires less compute time 
than the Crank-Nicholson scheme.

\section{\label{sec:applications_multidim}
         Application to PDEs in Multiple Space Dimensions}
OTS-NIDC can be used to boost the order of accuracy for finite difference 
schemes in any number of space dimensions.  However, several important issues 
arise for problems in more than one space dimension.  These issues are 
primarily related to the presence of cross-terms in the discretization error, 
multiple grid spacing parameters,  and the domain geometry.  In this section, 
we show how OTS-NIDC selection easily overcomes these issues by using it to 
construct high-order finite difference schemes for the 2D advection equation 
and the 2D diffusion equation.  
It is straightforward to generalize the approach presented in this section to 
more space dimensions and to construct methods for semilinear PDEs by 
incorporating ideas discussed in 
Sections~\ref{sec:OTS-NIDC} and~\ref{sec:applications_1d}.

\subsection{\label{sec:ots_multidim_crossterms}
            Discretization of Spatial Derivatives}
A significant distinction between the error for finite difference schemes in 
one space dimension and multiple space dimensions is the presence of 
cross-terms in the discretization error.  When formally analyzing the accuracy
of a finite difference scheme, it is usually safe to ignore the details of 
these terms because they can be absorbed within big-$O$ notation.
However, precise knowledge of the cross-terms is critical when attempting to 
use OTS-NIDC.  Without this knowledge, it is impossible to 
ensure cancellation of leading-order errors.  

Recall that OTS-NIDC requires the leading-order temporal errors to be directly 
related to the leading-order spatial errors via the PDE.  This requirement can 
be satisfied in one of two ways:
\begin{itemize}
\item by carefully designing the spatial discretization so that the
      leading-order spatial errors are directly related to the spatial 
      derivative operators in the PDE; or
\item by explicitly adding defect correction terms to cancel out any 
      cross-terms in the discretization error that cannot be eliminated 
      by using the PDE.
\end{itemize}
These two approaches are actually equivalent because the correction terms 
can be combined with the finite difference stencils in the original scheme 
and be viewed as part of an alternate spatial discretization scheme.

\subsection{Linear Advection Equation}
For our first application of OTS-NIDC in multiple space dimensions, we 
consider the two-dimensional advection equation
\beq
  \frac{\partial u}{\partial t} + \vec{A} \cdot \nabla u = 0,
  \label{eq:advection_eqn_2d}
\eeq
where $\vec{A} = (A_x, A_y)$ is the advection velocity.  
The most interesting and important aspect of this problem is that elimination 
of the truncation error requires optimal selection of the grid spacing ratio 
as well as the time step.  That is, the more general idea of numerical 
parameter optimization is needed in order to boost the accuracy of the scheme.

For convenience, let us assume that $A_x, A_y > 0$.  The simple forward Euler 
scheme for (\ref{eq:advection_eqn_2d}) with first-order upwind discretizations 
of the advection terms is 
\beq
  u^{n+1}_{i,j} = u^{n}_{i,j}
  - A_x \dt \left( \frac{u^{n}_{i,j} - u^{n}_{i-1,j}}{\dx} \right)
  - A_y \dt \left( \frac{u^{n}_{i,j} - u^{n}_{i,j-1}}{\dy} \right)
  \label{eq:advection_eqn_2d_FD_scheme}
\eeq
This method is formally first-order in both space and time with a 
stability constraint of $\dt < \left( A_x/\dx + A_y/\dy \right)^{-1}$.

To compute the optimal time step and grid spacing ratio, we begin by analyzing 
the discretization error 
for~(\ref{eq:advection_eqn_2d_FD_scheme}).  As in our analysis of the finite
difference scheme for the one-dimensional advection equation, we retain all of 
the terms in the discretization error.  Using Taylor series expansions and
the PDE (\ref{eq:advection_eqn_2d}), we find that the true solution $\tu$ 
satisfies
\bea
  \tu^{n+1}_{i,j} &=& \tu^{n}_{i,j }
    + \sum_{k=1}^\infty \frac{\left( -\dt \right)^k}{k!} 
         \left( A_x \partial_x
              + A_y \partial_y \right)^k \tu^{n}.
  \label{eq:advection_eqn_2d_time_err}
  \\
  \frac{\tu^{n}_{i,j} - \tu^{n}_{i-1,j}}{\dx} &=& 
  \frac{\partial \tu}{\px} 
  - \frac{1}{\dx} \sum_{k=2}^\infty \frac{\left( -\dx \right)^k}{k!} 
       \frac{\partial^k \tu^n}{\px^k} 
  \label{eq:advection_eqn_2d_space_err_x}
  \\
  \frac{\tu^{n}_{i,j} - \tu^{n}_{i,j-1}}{\dy} &=&
  \frac{\partial \tu}{\py} 
  - \frac{1}{\dy} \sum_{k=2}^\infty \frac{\left( -\dy \right)^k}{k!} 
       \frac{\partial^k \tu^n}{\py^k} 
  \label{eq:advection_eqn_2d_space_err_y}.
\eea
Therefore, the local truncation error is given by:
\bea
  \tau^{n}_{i,j} &=&
      \frac{A_x \dt}{2} \frac{\partial^2 \tu}{\px^2}
      \left( \dx - A_x \dt \right)
    + \frac{A_y \dt}{2} \frac{\partial^2 \tu}{\py^2}
      \left( \dy - A_y \dt \right)
    - A_x A_y \dt^2 \frac{\partial^2 \tu}{\px \py}
   \nonumber \\
   &+& \frac{A_x\dt}{\dx} \sum_{k=3}^\infty \frac{\left( -\dx \right)^k}{k!} 
       \frac{\partial^k \tu^n}{\px^k} 
   + \frac{A_y\dt}{\dy} \sum_{k=3}^\infty \frac{\left( -\dy \right)^k}{k!} 
       \frac{\partial^k \tu^n}{\py^k} 
   \nonumber \\
   &-& \sum_{k=3}^\infty \frac{\left( -\dt \right)^k}{k!} 
       \left( A_x \partial_x + A_y \partial_y
              \right)^k \tu^{n} 
  \label{eq:advection_eqn_2d_trunc_err}
\eea
From this expression, we see that choosing the ratio of the grid spacings to
be $\dx/\dy = A_x/A_y$ and setting the time step to be 
$\dt_{opt} = \dx/A_x = \dy/A_y$ eliminates two of the three leading-order 
terms in the discretization error\footnote{Note that the use of different 
values for $\dx$ and $\dy$ is consistent with physical intuition.  When 
computing the upwind derivative in a given coordinate direction, the upstream 
value we use should be at a distance that is proportional to the flow speed in 
that direction.  Otherwise, we will be using upstream data that is too far or 
too near in one of the coordinate directions.}.
The last leading-order term can be eliminated by adding the defect correction 
term
\bea
   A_x A_y \dt^2 \frac{\partial^2 \tu}{\px \py}.
  \label{eq:advection_eqn_2d_corr_term}
\eea
These choices enable us to boost the order of accuracy for the finite 
difference scheme (\ref{eq:advection_eqn_2d_FD_scheme}) from first- to 
second-order.

\begin{figure}[tb]
\begin{center}
\scalebox{0.35}{\includegraphics{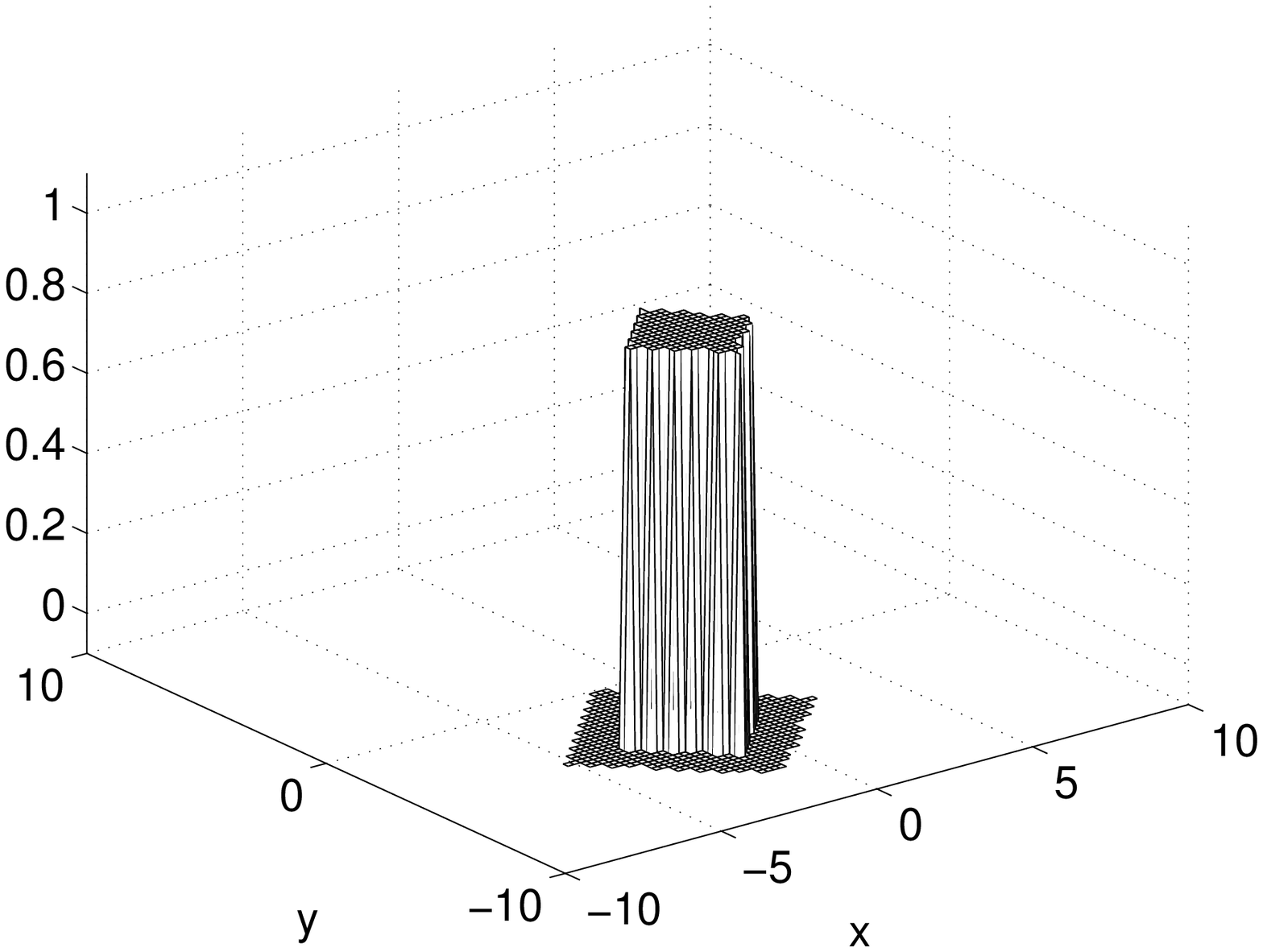}} 
\scalebox{0.35}{\includegraphics{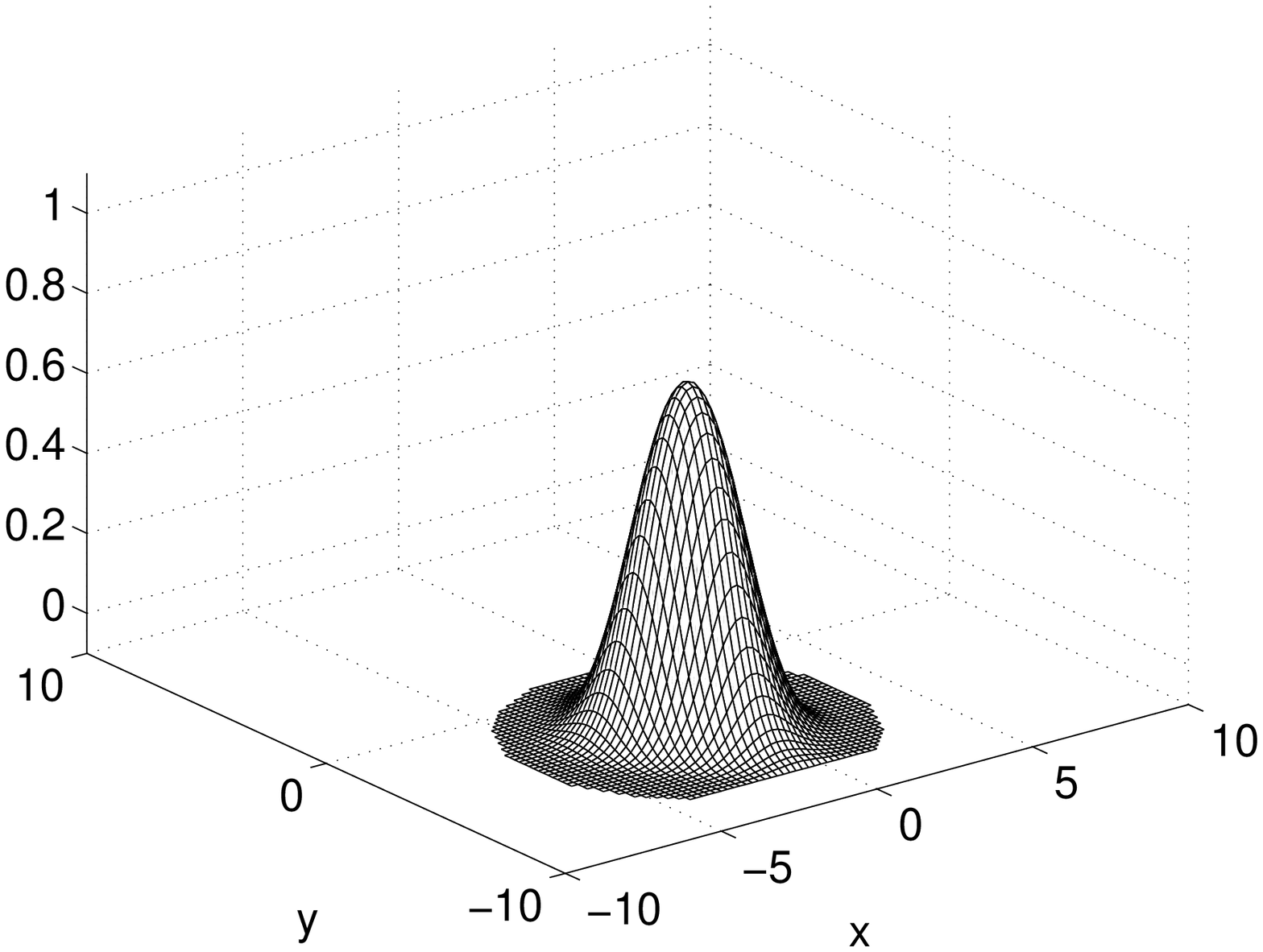}} 
\caption{Comparison of numerical solutions for the 2D advection
equation~(\ref{eq:advection_eqn_2d}) using forward Euler time stepping 
with OTS-NIDC (left) and forward Euler time 
stepping with suboptimal time step (right).  The latter scheme uses a 
time step $\Delta t = \Delta x / 2 (|A_x|+|A_y|)$.  
Both figures show the numerical solution at $t = 3$ computed on the domain 
$-10 < x,y < 10$ using 100 grid points in each direction.
The advection velocity $(A_x, A_y)$ is set to $(-1, -2)$.  The initial 
conditions are taken to be a 2D square function centered at the
origin: $u(x,y,0) = 1$ if $|x| + |y| \le 2$ and 
$u(x,y,0) = 0$ if $|x| + |y| > 2$. 
For the time interval considered, $u$ is identically zero 
on the boundary.
}
\label{fig:advection_eqn_2d_soln}
\end{center}
\end{figure}

With a careful choice of discretization for the correction term, we can
boost the order of accuracy even further.  Using the first-order upwind 
discretization for $\tu_{xy}$
\beq
   \frac{ \left(u_{i,j} - u_{i-1,j}\right) 
        - \left(u_{i,j-1} - u_{i-1,j-1}\right)}
        {\dx \dy},
  \label{eq:advection_eqn_2d_upwind_corr_approx}
\eeq
the truncation error~(\ref{eq:advection_eqn_2d_trunc_err}) becomes
\bea
  \tau^{n}_{i,j} &=&
      \frac{A_x \dt}{2} \frac{\partial^2 \tu}{\px^2}
      \left( \dx - A_x \dt \right)
    + \frac{A_y \dt}{2} \frac{\partial^2 \tu}{\py^2}
      \left( \dy - A_y \dt \right)
   \nonumber \\
   &+& \frac{A_x\dt}{\dx \dy} ( \dy - A_y \dt )
       \sum_{k=3}^\infty \frac{\left( -\dx \right)^k}{k!} 
              \frac{\partial^k \tu^n}{\px^k} 
   \nonumber \\
   &+& \frac{A_y\dt}{\dx \dy} ( \dx - A_x \dt )
       \sum_{k=3}^\infty \frac{\left( -\dy \right)^k}{k!} 
              \frac{\partial^k \tu^n}{\py^k} 
   \nonumber \\
   &+& \frac{A_x A_y\dt^2}{\dx \dy} 
       \sum_{k=3}^\infty \sum_{l=0}^{k}
              (-1)^k {k \choose l}
              \frac{\dx^l \dy^{k-l}}{k!} 
              \frac{\partial^k \tu^n}{\px^l \py^{k-l}} 
   \nonumber \\
   &-& \sum_{k=3}^\infty \frac{\left( -\dt \right)^k}{k!} 
       \left( A_x \partial_x + A_y \partial_y
              \right)^k \tu^{n} 
  \label{eq:advection_eqn_2d_trunc_err_mod},
\eea
which vanishes when the optimal time step and grid spacing ratio are used.
In other words, when the defect correction 
(\ref{eq:advection_eqn_2d_upwind_corr_approx}) is used 
in conjunction with the optimal time step and grid spacing ratio, we obtain
a finite difference scheme that introduces no truncation error at each time 
step.  Figure~\ref{fig:advection_eqn_2d_soln} demonstrates the amazing accuracy 
of this scheme over the first-order forward Euler 
scheme~(\ref{eq:advection_eqn_2d_FD_scheme}) without any modifications.
Notice how the sharp corners of the initial conditions are perfectly 
preserved when the optimal time step, optimal grid spacing ratio, and 
correction term are used.  In contrast, the sharp corners are
completely smoothed out in the solution computed using the standard 
first-order upwind scheme. 

It is interesting to push the analysis of the optimal finite difference 
scheme a little further and examine the complete finite difference scheme
(including defect correction term):
\bea
  u^{n+1}_{i,j} = u^{n}_{i,j}
  &-& A_x \dt \left( \frac{u^{n}_{i,j} - u^{n}_{i-1,j}}{\dx} \right)
  - A_y \dt \left( \frac{u^{n}_{i,j} - u^{n}_{i,j-1}}{\dy} \right)
  \nonumber \\
  &+& A_x A_y \dt^2 
        \left( \frac{u_{i,j} - u_{i-1,j} - u_{i,j-1} + u_{i-1,j-1}}
                    {\dx \dy} 
        \right).
  \label{eq:advection_eqn_2d_FD_scheme_complete}
\eea
Using the relationship between the optimal time step and the grid spacing
ratio, this scheme can be rewritten as
\bea
  u^{n+1}_{i,j} = u^{n}_{i,j}
  &-& \frac{A_x \dt}{2}
    \left( \frac{u^{n}_{i,j} - u^{n}_{i-1,j}}{\dx} 
         + \frac{u^{n}_{i,j-1} - u^{n}_{i-1,j-1}}{\dx} 
    \right)
  \nonumber \\
  &-& \frac{A_y \dt}{2}
    \left( \frac{u^{n}_{i,j} - u^{n}_{i,j-1}}{\dy} 
         + \frac{u^{n}_{i-1,j} - u^{n}_{i-1,j-1}}{\dy} 
    \right),
  \label{eq:advection_eqn_2d_FD_scheme_simplified}
\eea
which is essentially a forward Euler time integration scheme with 
first-order upwind discretizations of the spatial derivatives that are 
computed using all four upwind neighbors.  This formulation of 
(\ref{eq:advection_eqn_2d_FD_scheme_complete}) suggests that the highest
accuracy finite difference scheme for the 2D linear advection equation arises
by using a fully upwind discretization of the advection terms with the 
optimal time step and optimal ratio of grid spacings.

In higher space dimensions, it becomes tedious to carry out the analysis 
to derive the infinite order accurate finite difference scheme for the 
advection equation because many high-order correction terms must be retained.  
However, it is straightforward to generalize 
(\ref{eq:advection_eqn_2d_FD_scheme_simplified}) to $d$ space dimensions by 
using first-order discretizations for the spatial derivatives that are 
averages over the $2^{d-1}$ upwind discretizations in each coordinate 
direction.   As with the two-dimensional problem, the optimal grid spacings 
satisfy $\dx_k/A_{x_k} = \dt_{opt}$ for $1 \le k \le d$.

Like the unit CFL condition for the one-dimensional advection equation, the 
analogue of (\ref{eq:advection_eqn_2d_FD_scheme_simplified}) in any number
of space dimensions could also have been derived by examining how 
characteristic lines of the $d$-dimensional advection equation pass 
through grid points.  For the two-dimensional problem, this perspective 
becomes apparent when we fully simplify 
(\ref{eq:advection_eqn_2d_FD_scheme_simplified}) using the 
optimal time step and grid spacing ratio to obtain
$u^{n+1}_{i,j} = u^{n}_{i-1,j-1}$.

\subsection{Diffusion Equation \label{sec:diffusion_eqn_2d}}
As for PDEs in one space dimension, achieving infinite order accuracy 
through OTS-NIDC is unique to the advection equation.
To illustrate the typical accuracy boost gained by using OTS-NIDC,
we apply it to the 2D diffusion equation.  This example underscores 
the need to carefully choose the discretization for the spatial derivatives 
for problems in more than one space dimension and demonstrates the ease with 
which OTS-NIDC selection can be applied to problems on irregular domains.

In two space dimensions, the diffusion equation may be written as
\beq
  \frac{\partial u}{\pt} = D \nabla^2 u + f(x,y,t).
  \label{eq:diffusion_eqn_2d}
\eeq
The most common finite difference schemes for (\ref{eq:diffusion_eqn_2d}) 
employ a five-point, second-order central difference stencil for the 
Laplacian:
\beq
  L^{5pt} u = \frac{u_{i,j+1} + u_{i,j-1}
                   +u_{i+1,j} + u_{i-1,j} - 4u_{i,j}}{\dx^2},
  \label{eq:laplacian_2d_5pt_stencil}
\eeq
which assumes that the grid spacing in the $x$ and $y$ directions are
equal.  The forward Euler scheme using the five-point stencil is given by
\beq
  u^{n+1}_j = u^{n}_j 
  + \dt 
    \left( D \left[\frac{u^{n}_{i,j+1} + u^{n}_{i,j-1}
             +u^{n}_{i+1,j} + u^{n}_{i-1,j} - 4u^{n}_{i,j}}{\dx^2} \right]
         + f_{i,j}^n
    \right).
  \label{eq:diffusion_eqn_2d_5pt_scheme}
\eeq
This scheme is formally first-order in time and second-order in space.  
The stability constraint $\dt \le \dx^2/4D$ implies that this scheme is
$O(\dx^2)$ accurate overall.

While it is certainly possible to directly apply OTS-NIDC
to (\ref{eq:diffusion_eqn_2d_5pt_scheme}), we would be forced to use additional
defect correction terms to deal with the fact that the leading-order error for 
(\ref{eq:laplacian_2d_5pt_stencil}), 
$\dx^2 \left(\tu_{xxxx} + \tu_{yyyy}\right)/12$, cannot be directly related 
to the spatial derivative operator in the diffusion equation.  
The analysis is more straightforward if we construct our finite difference 
scheme using the nine-point stencil for the 
Laplacian~\cite{iserles_book,patra_2005}, which possesses a leading-order 
error directly related to the spatial derivative operator in the PDE 
(\ref{eq:diffusion_eqn_2d}):
\beq
  L^{9pt} u = \frac{1}{6 \dx^2} \left( 
      \begin{array}{l}
         \ \ \ u_{i+1,j+1} + u_{i+1,j-1}
            + u_{i-1,j+1} + u_{i-1,j-1} \\
         +\ 4 u_{i+1,j} + 4 u_{i-1,j}
           + 4 u_{i,j+1} + 4 u_{i,j-1}
           -20 u_{i,j} 
      \end{array}
    \right). 
  \label{eq:laplacian_2d_9pt_stencil}
\eeq
Observe that the nine-point stencil is simply the standard five-point 
stencil plus $\dx^2/6$ times a second-order, central difference 
approximation to $u_{xxyy}$.  Therefore, we can obtain the leading-order 
error in the nine-point stencil by combining the leading-order error for the 
five-point stencil with the extra contribution of 
$\dx^2 \tu_{xxyy}/6$
\bea
L^{9pt} \tu &=& \nabla^2 \tu 
              + \frac{\dx^2}{12} 
                \left(\tu_{xxxx} + 2 \tu_{xxyy} + \tu_{yyyy}\right)
              + O(\dx^4)
          \nonumber \\
          &=& \nabla^2 \tu + \frac{\dx^2}{12} \nabla^4 \tu + O(\dx^4).
  \label{eq:laplacian_2d_9pt_stencil_error}
\eea
Like (\ref{eq:diffusion_eqn_2d_5pt_scheme}), the forward Euler scheme that 
uses the nine-point discretization of the Laplacian 
\beq
  u^{n+1}_{i,j} = u^{n}_{i,j}
  + \dt \left( D L^{9pt} u + f_{i,j}^n \right)
  \label{eq:diffusion_eqn_2d_9pt_scheme}
\eeq
is $O(\dx^2)$ accurate overall because it has a stability constraint 
$\dt \le 3\dx^2/8D = O(\dx^2)$.
It is interesting to note that use of a nine-point stencil for the Laplacian 
naturally arises when directly applying OTS-NIDC to the five-point scheme if 
the defect correction term is interpreted as a modification to the original 
spatial discretization.

To apply OTS-NIDC to (\ref{eq:diffusion_eqn_2d_9pt_scheme}), we proceed in the 
usual way by first analyzing the leading-order terms of the local truncation
error for the finite difference scheme:
\bea
  \nabla^4 \tu
      \left[ \frac{\dx^2}{12} - \frac{D \dt}{2}  \right] (D \dt)
  - \frac{\dt^2}{2} \left( D\nabla^2 f + \frac{\partial f}{\pt} \right)
  + O(\dt^3) + O(\dt \dx^4).
  \label{eq:diffusion_eqn_2d_trunc_err}
\eea
Therefore, the optimal time step is given by $\dto = \dx^2/6D$ and the 
defect correction term is $\dt^2 \left( D\nabla^2 f + f_t \right)$.
This choice of time step and correction term yields a scheme that is 
fourth-order accurate.  It is interesting to note that the optimal time step
for the forward Euler scheme for the diffusion equation is \emph{independent} 
of the number of spatial dimensions.

\begin{figure}[tb]
\begin{center}
\scalebox{0.33}{\includegraphics{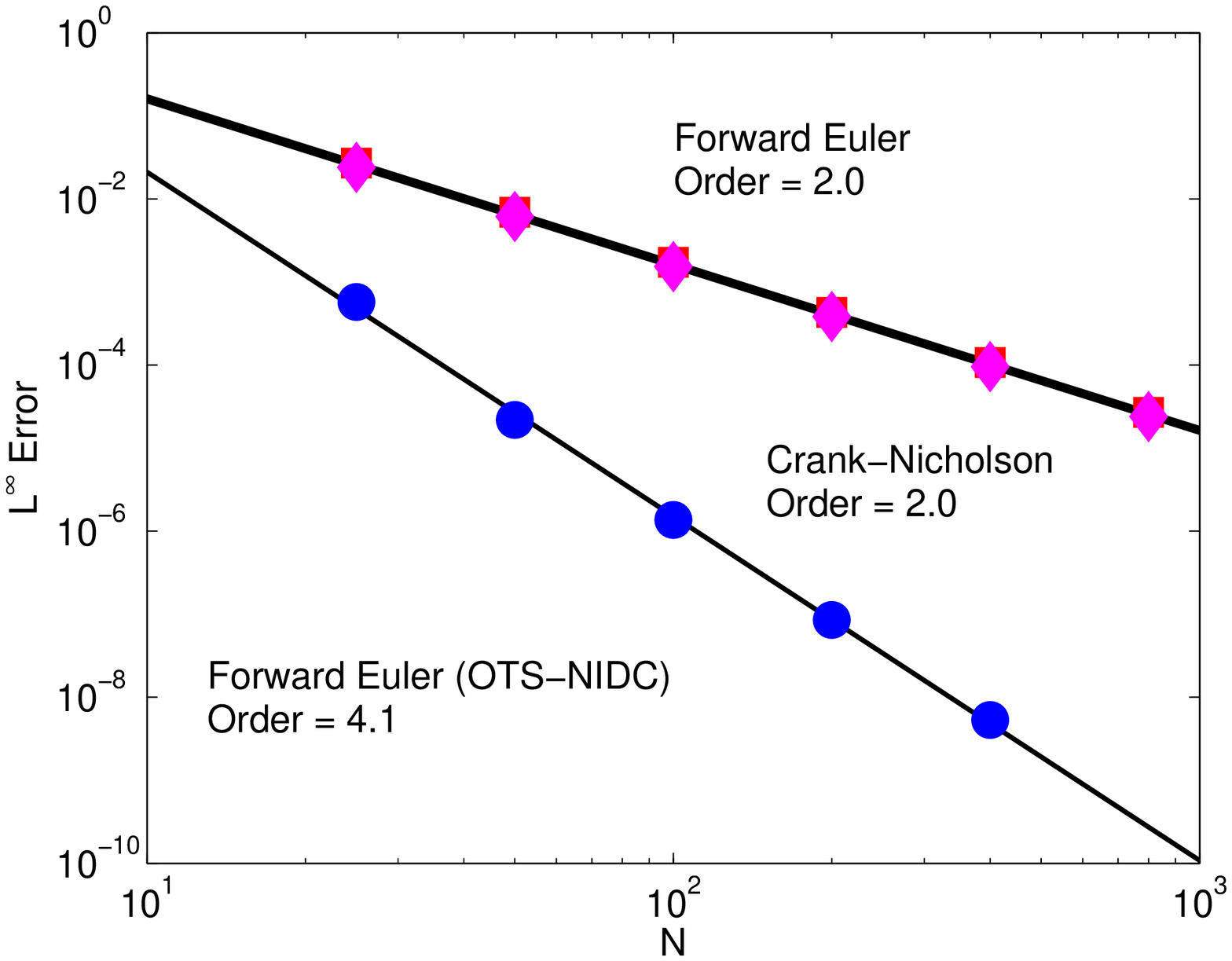}} 
\ \ 
\scalebox{0.33}{\includegraphics{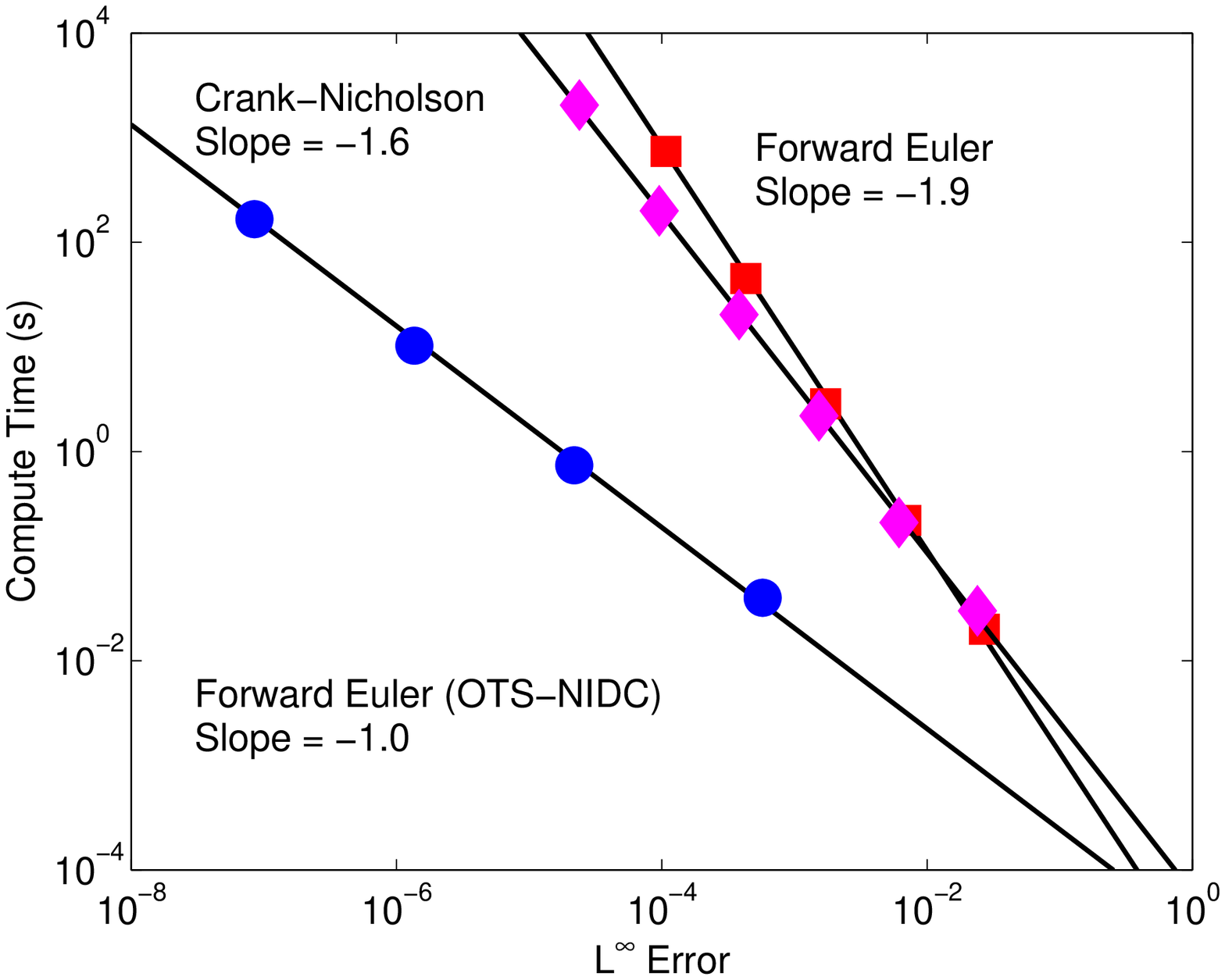}} 
\caption{$L^\infty$ error as a function of number of grid points in
each coordinate direction (left) and computation time as a function of 
$L^\infty$ error (right) for three finite difference schemes that solve the 
2D diffusion equation:
forward Euler with OTS-NIDC (circles), forward 
Euler with suboptimal time step (squares), and Crank-Nicholson (diamonds).  
The forward Euler schemes use a nine-point discretization of the Laplacian;
the Crank-Nicholson scheme uses the standard five-point discretization of
the Laplacian.  
Note that the errors for the forward Euler with suboptimal time step and 
Crank-Nicholson schemes lie almost directly on top of each other at the 
resolution of the figures.  
These results were obtained using MATLAB implementations of the 
finite difference schemes run on a 2.4 GHz MacBook Pro.
}
\label{fig:diffusion_eqn_2d_src_analysis}
\end{center}
\end{figure}
The improvement in the accuracy of the forward Euler scheme with OTS-NIDC
is shown in Figure~\ref{fig:diffusion_eqn_2d_src_analysis}.
Figure~\ref{fig:diffusion_eqn_2d_src_analysis} also shows that the 
scaling of the computational time with the error is in good agreement with the
theoretical estimates given in Table~\ref{tab:comp_perf_vs_dim}.  As with the
other examples, the forward Euler scheme with OTS-NIDC is much 
more computationally efficient than the two other schemes.

\subsubsection{Irregular Domains}
All of the examples we have considered so far share an important feature: they 
are all defined on simple, rectangular domains.  Any of several methods 
could have been used to obtain solutions with similar or even higher accuracy.
For example, a fourth-order accurate solution to the diffusion equation
(in any number of space dimensions) could be obtained by using a 
pseudospectral spatial 
discretization~\cite{trefethen_spectral_book,boyd_spectral_book} 
of the Laplacian and a fourth-order time integration scheme 
(\eg Runge-Kutta).  When the domain is irregularly shaped, high-order methods 
are harder to design and implement.  
In this section, we extend the OTS-NIDC finite difference scheme for the 2D 
diffusion equation on regular domains to irregular domains.  Other 
researchers have developed fourth-order accurate methods for the diffusion 
equation on irregular domains~\cite{gibou_2005,ito_2005}.  However, these 
earlier methods require a wider ghost cell region~\cite{gibou_2005}, use
implicit time integration~\cite{gibou_2005,ito_2005}, or involve a more 
complex calculation for the boundary conditions~\cite{ito_2005}.  As we shall 
see, it is possible to achieve fourth-order accuracy for the diffusion 
equation on irregular domains using only forward Euler integration and a 
compact, second-order stencil for the Laplacian.

% FIGURE SHOWING GHOST CELLS AT IRREGULAR BOUNDARY
% PLACED HERE TO SPREAD FIGURES OUT
\begin{figure}[tb]
\begin{center}
\scalebox{0.25}{\includegraphics{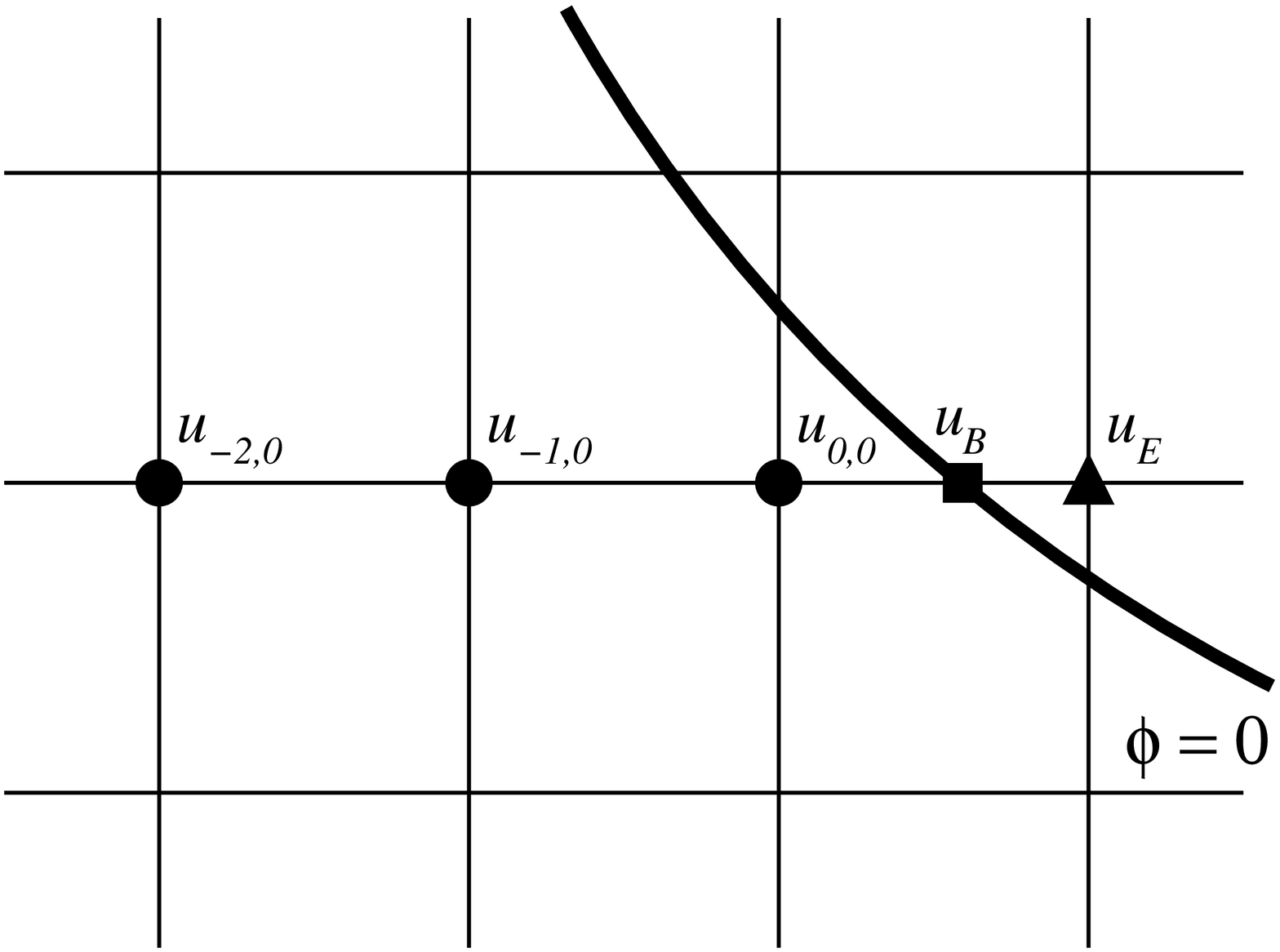}} 
\ \ \ \ \ \ \ \ \ \ \ \ \ 
\scalebox{0.25}{\includegraphics{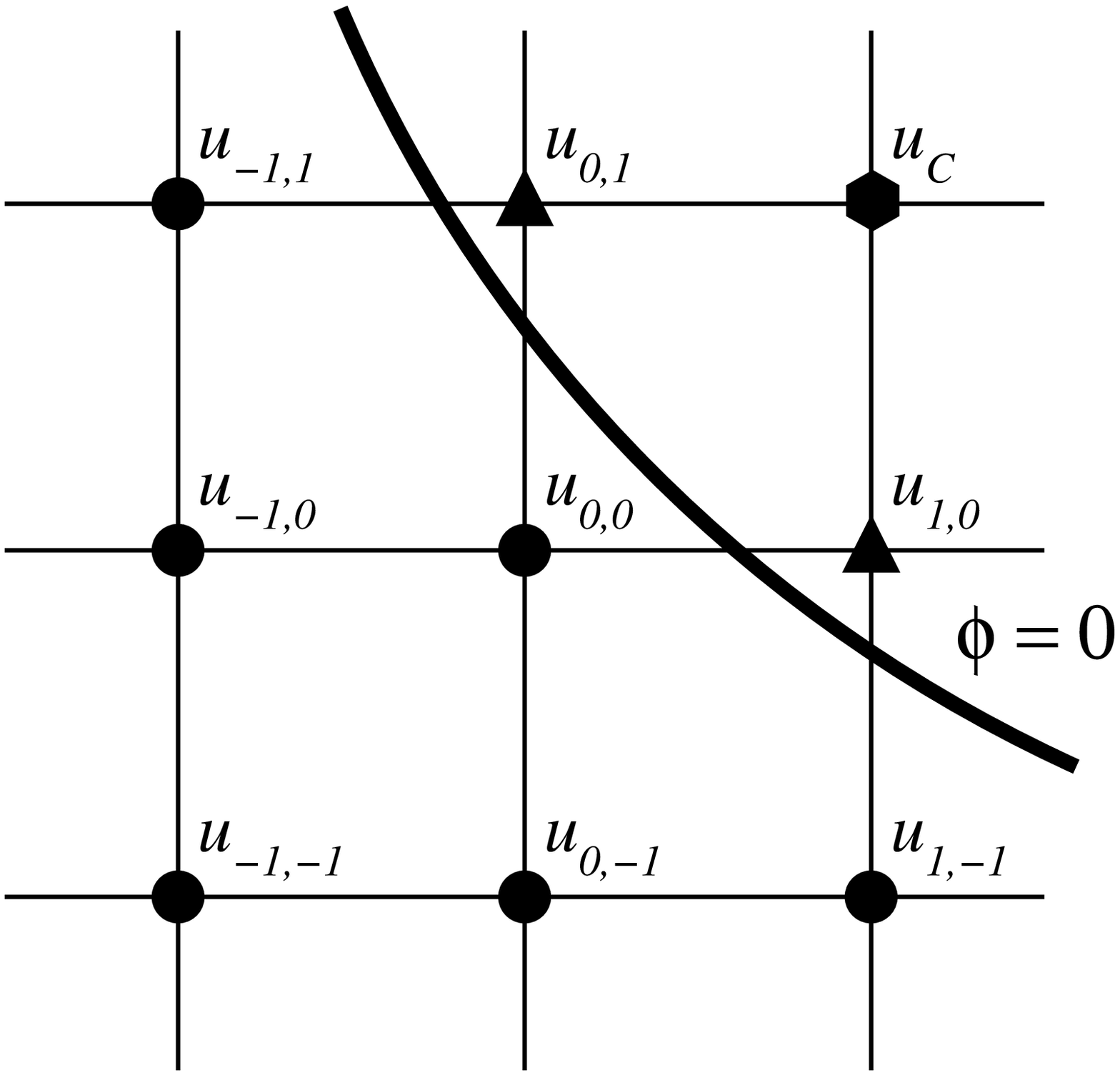}}
\caption{Illustrations of edge (left) and corner (right) ghost cells.
The edge ghost cell, $u_E$, is filled by using cubic extrapolation of 
the value on the boundary, $u_B$, and values at three interior grid points.  
The corner ghost cell, $u_C$, is filled by using a fourth-order accurate
Taylor series expansion centered at $u_{0,0}$, which uses the values of $u$ 
at the specified interior grid points and neighboring edge ghost cells to 
approximate the partial derivatives in the expansion. 
}
\label{fig:ghost_cells}
\end{center}
\end{figure}

To avoid the need for special stencils at grid points near the irregular 
boundary, we adopt the `ghost cell' approach for imposing boundary 
conditions~\cite{gibou_2005,ito_2005,fedkiw_1999,osher_fedkiw_book}.
The main challenge in using the ghost cell method is ensuring ghost cells 
are filled using a sufficiently high-order extrapolation scheme. 
Figure~\ref{fig:ghost_cells} shows examples of ghost cells that lie in the 
vicinity of an irregular boundary.  Notice that there are two types of ghost 
cells:  edge and corner.  Edge and corner ghost cells are distinguished by 
the relative position of the nearest interior grid point.  An edge ghost cell 
is linked to its nearest interior neighbor by an edge of the computational 
grid whereas a corner ghost cell and its nearest interior neighbor lie on the 
opposite corners of a grid cell.

To set the value of $u$ in each edge ghost cell, we use 1D cubic 
(\ie fourth-order accurate) extrapolation of the values of $u$ across the 
interface.  A cubic Lagrange interpolant is constructed for each edge ghost 
cell using a point on the domain boundary and three interior grid points that 
lie along the line parallel to the edge connecting the ghost cell to its 
nearest interior neighbor (see Figure~\ref{fig:ghost_cells})\footnote{If
a level set representation of the boundary~\cite{osher_fedkiw_book} is used, 
the location of the point $u_B$ can be calculated by inverting the linear 
interpolant of $\phi$ that passes through the ghost cell and its nearest
interior neighbor.  We found that using higher-order interpolants to 
locate boundary points was unnecessary for achieving high-order accuracy in 
the computed solution.}.  Because the quality of the interpolant deteriorates 
if the boundary point $u_B$ is too close to any of the interior grid points 
used to construct the Lagrange extrapolant, we follow \cite{gibou_2005} and 
choose the interior grid points so that the nearest interior grid point is 
sufficiently far from the boundary point.  Without this procedure, the errors 
introduced in the values of the ghost cells lead to instability of the 
numerical method.
Gibou~\etal~suggest shifting the interpolation points by one grid point 
when the distance between the boundary point and the nearest interior grid 
point is less than $\dx^2$~\cite{gibou_2005}.  For the numerical scheme 
presented in this article, this threshold was too low.  Better stability was 
obtained if an $O(\dx)$ threshold was used.

% FIGURE SHOWING SOLUTION OF 2D DIFFUSION EQUATION ON STARFISH DOMAIN
% PLACED HERE TO SPREAD FIGURES OUT
\begin{figure}[tb]
\begin{center}
\scalebox{0.35}{\includegraphics{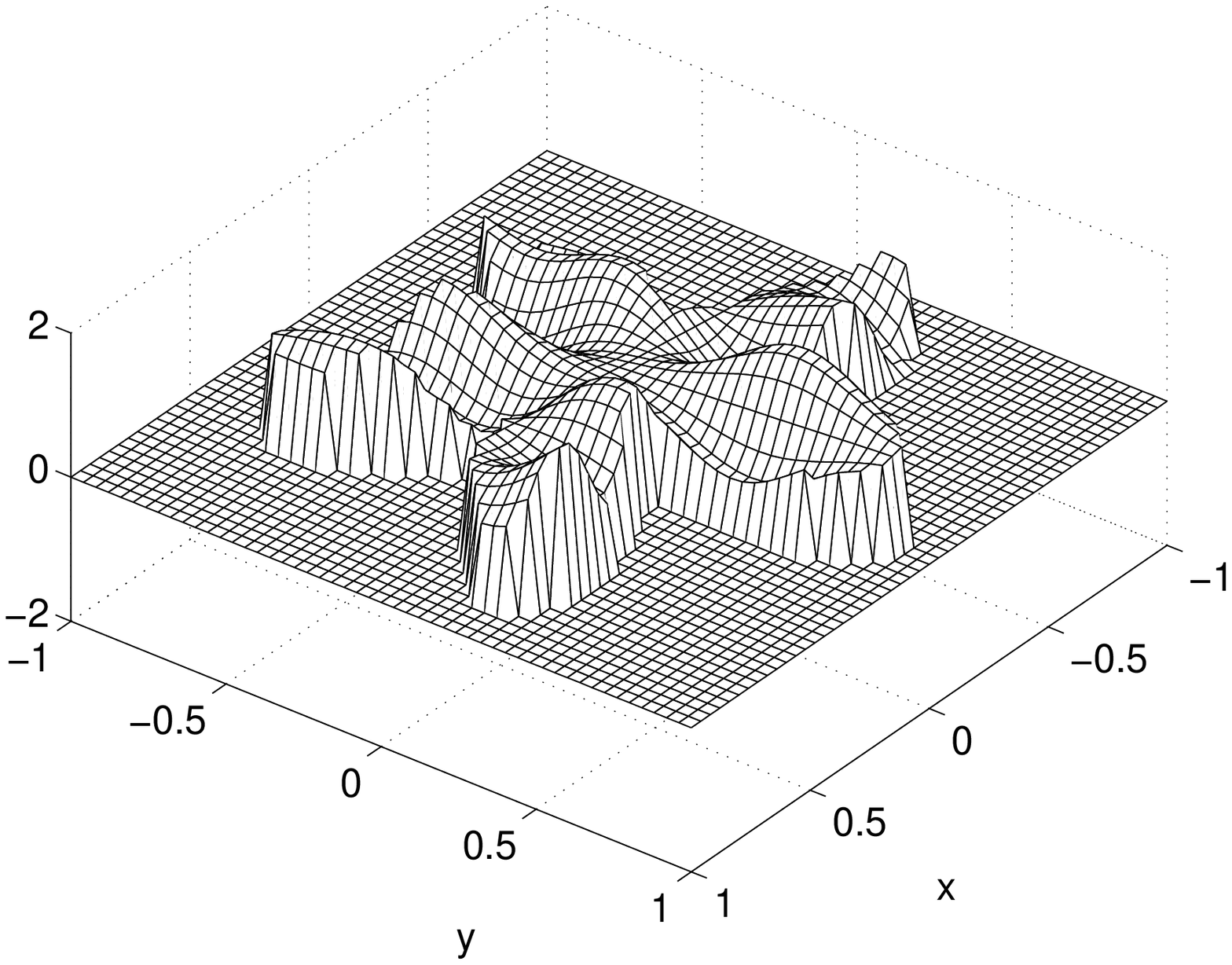}} 
\scalebox{0.35}{\includegraphics{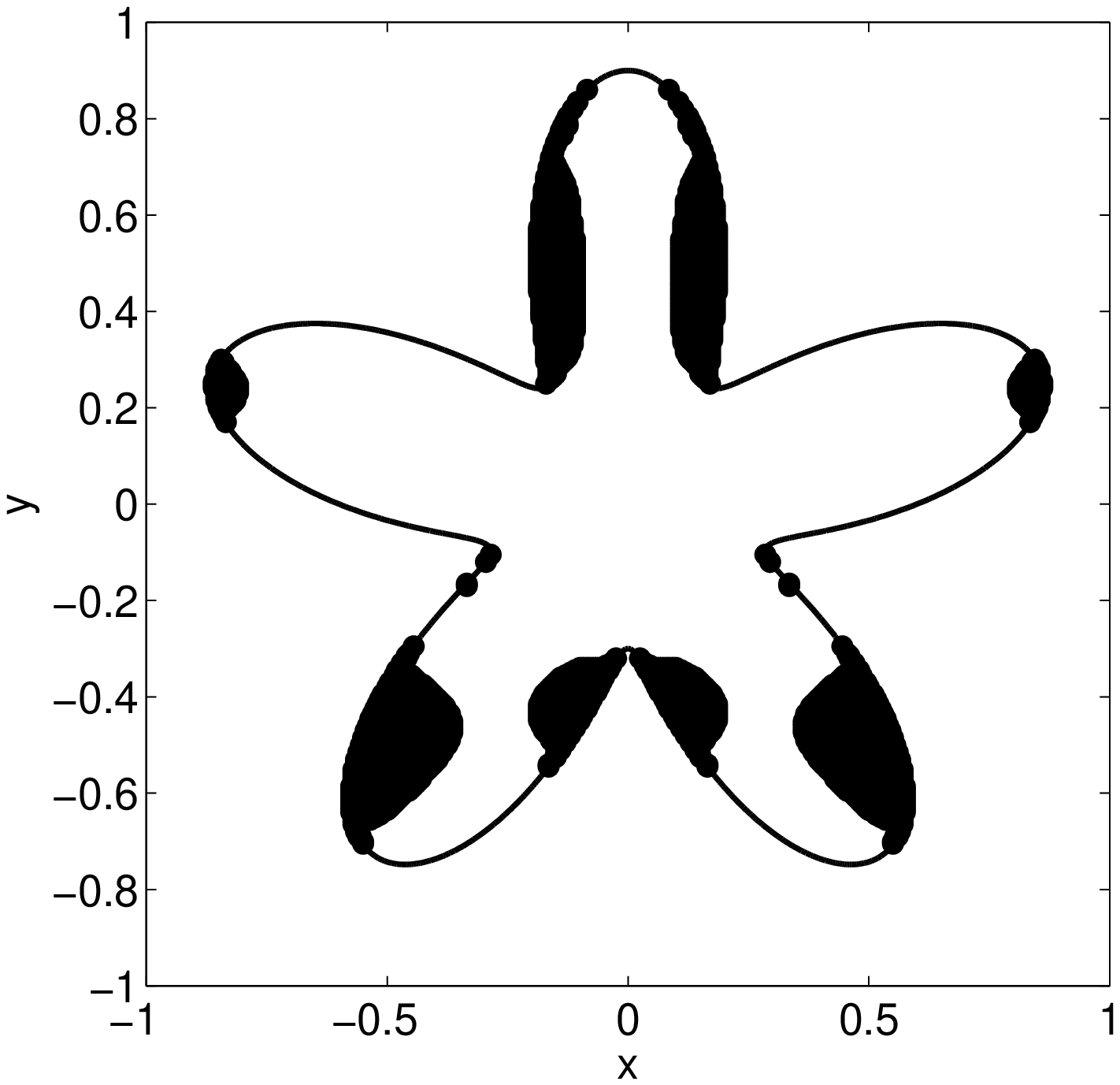}} 
\caption{Numerical solutions (left) and dominant error (right) for the 2D 
diffusion equation on a starfish-shaped domain.  The solution is computed 
using forward Euler time stepping with OTS-NIDC.  In the error plot, the dark 
regions represent points where the error in the solution is larger than 
$25$\% of the $L^\infty$ error of the solution.  
}
\label{fig:diffusion_eqn_2d_starfish_domain}
\end{center}
\end{figure}

To set the value of $u$ in a corner ghost cell, we use a fourth-order accurate
Taylor series expansion of $u$ centered at the nearest interior 
neighbor.  The partial derivatives in the Taylor series expansion are 
approximated using finite differences constructed from the values of $u$ 
at interior grid points and the neighboring edge ghost cells (see 
Figure~\ref{fig:ghost_cells}).  For the corner ghost cell shown 
in Figure~\ref{fig:ghost_cells}, the extrapolation stencil is given by
\bea
  u_C &=&  -4 u_{0,0} - u_{-1,-1} + 2 u_{1,0} + 2 u_{0,1}
      - u_{1,-1} - u_{-1,1} + 2 u_{0,-1} + 2 u_{-1,0}
\eea
where $u_C$ and $u_{i,j}$ are indicated in the figure.  The stencils for
corner ghost cells that are in other positions relative to the interior of 
the domain are obtained by rotations of Figure~\ref{fig:ghost_cells} and 
remapping of the indices in the above stencil.

\begin{figure}[tb]
\begin{center}
\scalebox{0.35}{
  \includegraphics{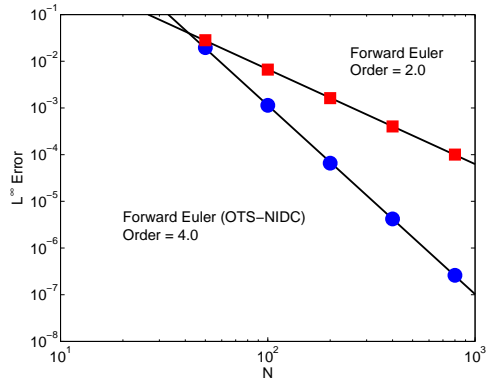}} 
\caption{$L^\infty$ error as a function of number of grid points in each
coordinate direction for two finite difference schemes that solve the 2D 
diffusion equation on an irregular domain: forward Euler with OTS-NIDC 
(circles) and forward Euler with suboptimal time step (squares).  
}
\label{fig:diffusion_eqn_2d_starfish_error}
\end{center}
\end{figure}

Using these methods for filling edge and corner ghost cells, it is 
straightforward to compute high-order accurate solutions of the 2D diffusion 
equation on irregular domains.  We simply use the OTS-NIDC forward Euler 
scheme designed for regular domains for grid points in the interior of the 
domain after filling the ghost cells with high-order accurate values.  
Figure~\ref{fig:diffusion_eqn_2d_starfish_domain} shows an example of 
a solution computed on a starfish shaped domain.  It also shows that the 
error in the solution is concentrated near the boundaries of the irregular
domain.  The improved accuracy of the OTS-NIDC solution compared to a simple 
forward Euler scheme with a suboptimal time step is illustrated in 
Figure~\ref{fig:diffusion_eqn_2d_starfish_error}. 

% FIGURE SHOWING SOLUTION OF 2D DIFFUSION EQUATION ON STARFISH DOMAIN
% PLACED HERE TO SPREAD FIGURES OUT

\subsubsection{Diffusion Equation in Higher Dimensions}
Using OTS-NIDC to boost the accuracy of finite difference
schemes for the diffusion equation in higher space dimensions is 
straightforward.  As for the 2D problem, the key step is identification of a 
finite difference stencil for the Laplacian that has an \emph{isotropic} 
discretization error.  For 3D problems, a catalog of several such stencils is 
available in~\cite{patra_2005}.   For higher dimensional problems, one simple
approach for deriving an isotropic finite difference stencil for the Laplacian
would be to start with the generalization of the standard five point stencil
for the 2D Laplacian and add finite difference approximations for the missing 
cross terms that are required to make the leading-order error isotropic 
(\ie proportional to the bilaplacian of the solution).  Because the missing 
terms are all of the form 
\beq
\frac{\dx^2}{6} \left(\frac{\partial^4 u}{\partial x_i^2 \partial x_j^2}\right),
\eeq
this approach is relatively easy to use.  To derive finite difference 
approximations with more general or complex properties, symbolic algebra 
software can be helpful~\cite{patra_2005,gupta_1998}.  

Once an acceptable finite difference approximation for the Laplacian has been 
chosen, the usual OTS analysis yields an optimal time step of 
$\dt_{opt} = \dx^2/6D$ for forward Euler time integration and a boost of the 
order of accuracy from two to four.  Above three dimensions, the optimal
time step is no longer stable, and the forward Euler scheme should be replaced
with the DuFort-Frankel scheme (see Section~\ref{sec:dufort_frankel}).

Problems on irregular domains can be handled using the ghost cell approach. 
While the derivation of high accuracy stencils for ghost cells becomes more
tedious, the general approach is the same as for the 2D problem.  The ghost 
cells can first be grouped into separate classes based on the number of 
grid cell edges that must be traversed to reach the nearest interior neighbor.
The values in the ghost cells can then be set in order starting with those 
directly connected to their nearest interior neighbor and ending with those 
whose nearest interior neighbor is on the opposite corner of a grid cell.
Here, too, symbolic algebra programs can be helpful when deriving methods for 
extrapolating interior and boundary values to ghost cells.

\section{\label{sec:summary} Summary} 
In this article, we have presented OTS-NIDC, a simple approach for 
constructing high-order finite difference methods for time dependent PDEs 
\emph{without} requiring the use of high-order stencils or high-order time 
integration schemes.  The major idea underlying OTS-NIDC is that 
spatial and temporal discretization errors can be simultaneously eliminated
by (1) using the PDE to relate terms in the discretization error, (2)
optimally selecting numerical parameters, and (3) adding a few defect 
correction terms in a non-iterative manner.  The boost in the order of 
accuracy achieved by using OTS-NIDC is well worth the extra effort required to 
derive the optimal time step and defect correction terms. 

Through many examples, we have demonstrated the utility of OTS-NIDC
to several types of finite difference schemes for a wide range of problems.  
We showed how OTS-NIDC can be used to very easily obtain high-order accurate 
solutions to many linear and semilinear PDEs in any number of space dimensions 
on both regular and irregular domains.  The examples illustrate the 
high-order of accuracy that can be achieved using only low-order 
discretizations for spatial derivatives and simple forward Euler time 
integration.  They also demonstrate the applicability of OTS-NIDC to more 
exotic finite difference schemes, such as the Kreiss-Petersson-Ystr\"om 
discretization of the second-order wave equation and the DuFort-Frankel scheme 
for the diffusion equation.

Optimal time step selection is an example of the more general technique
of optimally selecting numerical parameters to boost the accuracy of 
numerical methods.  Little research appears to have been done in this area,
but as OTS-NIDC selection demonstrates, optimal selection of numerical 
parameters can significantly impact the accuracy of numerical methods.  To 
realize the full potential of existing and novel numerical schemes, 
both optimal time step selection and optimal numerical parameter selection 
deserve further investigation and attention.

% ACKNOWLEDGMENTS
\section*{Acknowledgments}
The author gratefully acknowledges the support of Vitamin D, Inc.
and the Institute for High-Performance Computing (IHPC) in Singapore. 
The would like to thank P. Fok and J. Lambers for many enlightening 
discussions and P. Fok, M. Prodanovi\'c and A. Chiu for helpful suggestions 
on the manuscript.

\appendix
\section{Formal vs. Practical Accuracy
         \label{app:formal_vs_practical_accuracy} }
For time dependent PDEs, evaluating the order of accuracy for a numerical
method is subtle because there are formally two separate orders of accuracy 
to consider:  temporal and spatial.  While it is theoretically
interesting and important to understand how the error depends on both the 
grid spacing $\dx$ and the time step $\dt$, in practice, the accuracy 
of a numerical scheme is always controlled by only one of the two numerical
parameters.  

Even though spatial and temporal orders of accuracy for a numerical method
are formally separate, one will always dominate for a given choice of 
$\dx$ and $\dt$.  For example, when solving the diffusion equation 
using the backward Euler method for time integration with a second-order 
central difference stencil for the Laplacian, formal analysis shows that 
the global error is $O(\dx^2) + O(\dt)$.  Since there are no 
stability constraints on the numerical parameters, the time
step and grid spacing are free to vary independently.  In this situation, the 
practical error for the method depends on the relative sizes of the time step 
and grid spacing.  When $\dt \gg \dx^2$, the practical error is 
$O(\dt)$ which means that the error in the numerical solution is 
primarily controlled by the time step.  Similarly, when 
$\dt \ll \dx^2$, the practical error is $O(\dx^2)$ so that 
the error is controlled by the grid spacing.  Finally, when 
$\dt  = O(\dx^2)$, the practical error is 
$O(\dt) = O(\dx^2)$.  In all cases, the practical error is 
primarily controlled by one of the two numerical parameters, and varying
the subdominant parameter while holding the dominant parameter fixed does 
not significantly affect the error.
 
When there are constraints on the numerical parameters, there is less freedom 
in choosing the controlling parameter.  For instance, if we solve the 
diffusion equation using a forward Euler time integration scheme with a 
second-order central difference stencil for the Laplacian, stability
considerations require that we choose $\dt = O(\dx^2)$.  
Combining this stability constraint with the formal error for the scheme
shows that the practical error is 
$O(\dx^2) + O(\dt) = O(\dx^2)$.  Therefore, the accuracy of the 
method is completely controlled by the grid spacing; the temporal error can 
never dominate the spatial error.

\section{Importance of High-Accuracy for First Time Step for the KPY Scheme
         \label{appendix:KPY_analysis} }
The need for higher-order accuracy when taking the first time step of the KPY 
scheme for the second-order wave equation can be understood by solving the 
difference equation 
\bea
  \he^{n+1} - 2 \he^{n} + \he^{n-1} = \dt^2 \lambda \he^{n}
  \label{eq:error_eqn_normal_mode},
\eea
for the normal modes of the error, where $\he$ is the coefficient of an 
arbitrary normal mode of the spatial operator for the error 
$e \equiv u - \tu$.  Using standard methods for solving linear difference
equations~\cite{carrier_pearson_book} yields the solution
\bea
  \he^n = \he^1 \frac{\kappa_+^n - \kappa_-^n}{\kappa_+ - \kappa_-}
        + \he^0 \kappa_+ \kappa_-
          \frac{\kappa_+^{n-1} - \kappa_-^{n-1}}{\kappa_+ - \kappa_-},
  \label{eq:error_eqn_normal_mode_soln}
\eea
where $\kappa_\pm$ are the roots of the characteristic equation for
(\ref{eq:error_eqn_normal_mode})~\cite{kreiss2002}
\beq
  \kappa_\pm = 1 + \frac{1}{2} \lambda \dt^2
             \pm \dt \sqrt{\lambda + \frac{\lambda^2 \dt^2}{4}}.
\eeq
Notice that the denominator of both terms in
(\ref{eq:error_eqn_normal_mode_soln}) is $O(\dt)$, which implies that the
error in the initial conditions is degraded by one temporal order of accuracy
by the KPY scheme.  Therefore, KPY without OTS-NIDC requires the first time 
step must be at least third-order accurate; KPY with OTS-NIDC requires the 
first time step must be at least fifth-order accurate.

%*********************************BIBLIOGRAPHY**********************************
%\bibliography{OTSPDE}

\end{document}